\documentclass[11pt,leqno]{article}
\usepackage{a4wide,amsmath}
\usepackage{amsfonts,amsbsy,amscd}
\usepackage{amssymb}
\usepackage{mathrsfs}
\usepackage{stmaryrd}
\usepackage{marvosym}
\usepackage{mathcomp}
\usepackage{color}
\usepackage{tcolorbox}
\usepackage{tikz}
\usepackage[colorlinks=true]{hyperref}
\newtheorem{theorem}{Theorem}[section]
\newtheorem{corollary}[theorem]{Corollary}
\newtheorem{lemma}[theorem]{Lemma}

\newtheorem{remark}[theorem]{Remark}
\newtheorem{example}[theorem]{Example}

\newtheorem{hypo}{Hypothesis}
\newtheorem{hypos}{Hypotheses}
\numberwithin{equation}{section}

\newcommand{\vertiii}[1]%
{{\left\vert\kern-0.25ex\left\vert\kern-0.25ex\left\vert #1
  \right\vert\kern-0.25ex\right\vert\kern-0.25ex\right\vert}}
\newcommand{\RR}{\mathbb{R}}
\newcommand{\CC}{\mathbb{C}}

\newcommand{\ee}{\mathrm{e}}
\newcommand{\ii}{\mathrm{i}}
\renewcommand{\colon}{:\,}

\newcommand{\eqdef}{\stackrel{{\rm {def}}}{=}}   % Composed symbols

\newcommand{\Square}{$\sqcap$\hskip -1.5ex $\sqcup$}
\newcommand{\Blacksquare}{\vrule height 1.7ex width 1.7ex depth 0.2ex }
\newcommand{\proof}{{\em Proof. }}
\newcommand{\qed}{$\;$\Blacksquare}

\definecolor{violet}{rgb}{0.5,0,0.5}
\definecolor{orange}{cmyk}{0,0.3,0.7,0}

\title{\null\vspace*{-2.0cm}
      Monotone methods in counterparty risk models\\
      with non\-linear Black\--Scholes\--type equations%
\footnote{A part of this research was performed while
          the second author (P.T.) was a visiting professor at
          Toulouse School of Economics, I.M.T.,
          Universit\'e de Toulouse -- Capitole,
          Toulouse, France.\vspace{2.0mm}}
\vspace*{0.5cm}
}
 
\author{%
        B\'en\'edicte Alziary%
\thanks{{\it e-mail:} {\tt benedicte.alziary@ut-capitole.fr},$\;$}
\vspace*{0.3cm}
\\
        Toulouse School of Economics, I.M.T.,
        Universit\'e de Toulouse -- Capitole \\
        1, Esplanade de l'Universit\'e,
        F--31000 Toulouse Cedex, France \\
\vspace*{0.3cm}
\\
        and \\
\vspace*{0.3cm}
\and
        Peter Tak\'a\v{c}%
\thanks{{\it e-mail:} {\tt peter.takac@uni-rostock.de}.$\;$}
\vspace*{0.3cm}
\\
        Universit\"at Rostock,
        Institut f\"ur Mathematik \\
        Ulmenstra{\ss}e~69, Haus~3,
        D--18057 Rostock, Germany
\vspace*{0.5cm}
}
\date\today

\begin{document}
\baselineskip=16pt plus 1pt minus 1pt

\maketitle
%\vspace{3cm}
\baselineskip=14pt plus 1pt minus 1pt
%\tableofcontents
\baselineskip=16pt plus 1pt minus 1pt
 
%\newpage

\begin{abstract}
A non\-linear Black\--Scholes\--type equation
is studied within {\em\bfseries counterparty risk models\/}.
The classical hypothesis on the uniform Lipschitz\--continuity
of the non\-linear reaction function allows for
an equivalent transformation of the semi\-linear Black\--Scholes equation
into a standard parabolic problem with
a monotone non\-linear reaction function and
an inhomogeneous linear diffusion equation.
This setting allows us to construct a scheme of monotone,
increasing or decreasing, iterations that converge monotonically
to the true solution.
As typically any numerical solution of this problem uses
most computational power for computing an approximate solution
to the inhomogeneous linear diffusion equation,
we discuss also this question and suggest several solution methods,
including those based on Monte Carlo and finite differences/elements.
\end{abstract}

\par\vfill
\vspace*{0.5cm}
\noindent
\begin{tabular}{lll}
{\bf 2020 Mathematics Subject Classification:}
& Primary   & 35A16, 91G40;\\
& Secondary & 35K58, 91G60;\\
\end{tabular}

\par\vspace*{0.5cm}
\noindent
\begin{tabular}{ll}
{\bf Key words:}
& non\-linear Black\--Scholes equation; counterparty risk models;\\
& semi\-linear parabolic Cauchy problem; monotone methods; \\
& non\-linear integral equation; fixed point by monotone iterations \\
\end{tabular}
 
%%%%%%%%%%%%%%%%%%%%%%%%%%%%%%%%%%%%%%%%%%%%%%%%%%%%%%%%%%%%%%%%%%%%%%%%
%\baselineskip=14pt plus 1pt minus 1pt
\baselineskip=16pt plus 1pt minus 1pt
\parskip=2mm plus .5mm minus .5mm
%%%%%%%%%%%%%%%%%%%%%%%%%%%%%%%%%%%%%%%%%%%%%%%%%%%%%%%%%%%%%%%%%%%%%%%%
\newpage
 
%%%%%%%%%%%%%%%%%%%%%%%%%%%%%%%%%%%%%%%%%%%%%%%%%%%%%%%%%%%%%%%%%%%%%%%
%%%%%    INTRODUCTION    %%%%%%%%%%%%%%%%%%%%%%%%%%%%%%%%%%%%%%%%%%%%%%
%%%%%%%%%%%%%%%%%%%%%%%%%%%%%%%%%%%%%%%%%%%%%%%%%%%%%%%%%%%%%%%%%%%%%%%

\section{Introduction}
\label{s:Intro}

Risk phenomena and their management have been an important topic
of investigation since the financial crisis of $2007$ -- $2008$.
In this article we focus our attention on
{\sl\bfseries counterparty risk models\/} for options with risky values
$\hat{V}(S,t)\in \mathbb{R}$
modelled by a non\-linear Black\--Scholes\--type equation:
\begin{align}
\label{e:BS_M=hatV}
&
\begin{aligned}
  \frac{\partial\hat{V}}{\partial t} + \mathcal{A}_t\hat{V}
  - (r + \lambda_B + \lambda_C) \hat{V}
  = F(\hat{V}(S,t); S,t)
\\
    \quad\mbox{ for }\, (S,t)\in (0,\infty)\times (0,T) \,;
\end{aligned}
\\
\label{e:BS_M=hatV,T}
& \hat{V}(S,T) = h(S)
    \quad\mbox{ for }\, S\in (0,\infty) \,.
\end{align}
The non\-linearity,
$F(\,\cdot\, ; S,t)\colon \RR\to \RR$,
with $\RR = (-\infty, +\infty)$ standing for the real line, is given by
\begin{align}
\label{e:F(hatV)}
&
\begin{aligned}
  F(M; S,t)\eqdef 
& {}= (R_B\lambda_B + \lambda_C)\, M^{-}
    - (\lambda_B + R_C\lambda_C)\, M^{+} + s_F\, M^{+}
\\
&   \quad\mbox{ for $M\in \RR$ and }\, (S,t)\in (0,\infty)\times (0,T) \,,
\end{aligned}
\end{align}
where we use the usual abbreviation
$x^{+}\eqdef \max\{ x,\,0\}$ and $x^{-}\eqdef \max\{ -x,\,0\}$
for $x\in \RR$.
Hence, $x = x^{+} - x^{-}$ and $|x| = x^{+} + x^{-}$.
These kinds of non\-linearities, often called
{\em ``jumping non\-linearities\/}, have a long tradition
in Mathematical Modelling.

The parabolic partial differential equation \eqref{e:BS_M=hatV}
(PDE, for short)
corresponds to the case when the non\-linearity
\begin{math}
  F(\,\cdot\,; S,t)\colon M\mapsto F(M; S,t)\colon \RR\to \RR
\end{math}
on the right\--hand side in eq.~\eqref{e:BS_M=hatV}
is taken with the {\it\bfseries mark\--to\--market\/} value
$M = \hat{V}(S,t)$.
This case corresponds to a derivative contract $\hat{V}$ on an asset
(stock)
$S\in (0,\infty)$ between a {\sl\bfseries seller\/}~\textbf{B}
and a {\sl\bfseries counterparty\/}~\textbf{C} that
{\sl\bfseries may both default\/}.
The asset price $S$ is not affected by
a default of either \textbf{B} or \textbf{C},
and is assumed to follow the Markov process with
the (time\--dependent) generator (the Black\--Scholes operator)
$\mathcal{A}_t$ defined by
\begin{equation}
\label{e:BS_gener}
\begin{aligned}
  (\mathcal{A}_t V)(S,t)\eqdef
  \frac{1}{2}\, [\sigma(t)]^2 S^2\, \frac{\partial^2 V}{\partial S^2}
  + [ q_S(t) - \gamma_S(t) ] S\, \frac{\partial V}{\partial S}
\\
    \quad\mbox{ for }\;
  V\colon (0,\infty)\times (0,T)\to \RR\colon
  (S,t)\mapsto V(S,t) \,.
\end{aligned}
\end{equation}
As usual, we take the volatility, $\sigma$, to be a positive constant,
$\sigma\in (0,\infty)$.
The value of $\gamma_S(t)$ reflects the
{\it\bfseries rate of dividend income\/}
and the value of $q_S(t)$ is the
{\it\bfseries net share position finan\-cing cost\/}
which depends on
the {\it\bfseries risk\--free rate\/} $r(t)$ and
the {\it\bfseries repo\--rate\/} of $S(t)$.
``Typical'' values for the terminal condition \eqref{e:BS_M=hatV,T}
are $h(S)\equiv V_T(S)$ where
$V_T(S) = (S-K)^{+} = ( \ee^X - K )^{+}$ for $X = \log~S\in \RR$
(in case of the {\em European call option\/})
and
$V_T(S) = (S-K)^{-} = (K-S)^{+} = ( K - \ee^X )^{+}$
(for the {\em European put option\/}).

A frequently used {\em alternative\/} to our choice
$M = \hat{V}(S,t)$ of the {\it\bfseries mark\--to\--market\/} value $M$
in the non\-linearity $F(M; S,t)$
on the right\--hand side in eq.~\eqref{e:BS_M=hatV}
is $M = V(S,t)$ where $V$ denotes the same derivative between two parties
that {\sl\bfseries cannot default\/}; see e.g.\
{\sc F.\ Baustian}, {\sc M.\ Fencl}, {\sc J.\ Posp\'{\i}\v{s}il}, and
{\sc V.\ {\v{S}}v\'{\i}gler} \cite[Sect.~2]{Fencl-Pospis-22}
for numerical treatment.
This {\em\bfseries risk\--free value\/}, $V$, satisfies
the classical (linear) Black\--Scholes PDE
(partial differential equation)
with the prescribed terminal value $V(S,T) = h(S)$
for $S\in (0,\infty)$ at {\it\bfseries maturity time\/} $t=T$.
Inserting this known value $F(V(S,t); S,t)$ in eq.~\eqref{e:BS_M=hatV}
in place of $F(\hat{V}(S,t); S,t)$,
we thus obtain an inhomogeneous linear equation for
another (new) value of $\hat{V}(S,t)$.
We refer to the works by
{\sc C.\ Burgard} and {\sc M.\ Kjaer}
\cite{Burgard-Kjaer},
\cite[Section~3]{Burgard-Kjaer-PDE}, and
\cite{Burgard-Kjaer-PDE-add} for details concerning modelling and to
{\sc I.\ Arregui\/}, {\sc B.\ Salvador\/}, and {\sc C.\ V{\'a}zquez\/}
\cite{Arregui-Salvad} for numerical results.
We warn the reader that Refs.\
\cite{Arregui-Salvad} and
\cite{Burgard-Kjaer, Burgard-Kjaer-PDE, Burgard-Kjaer-PDE-add}
use the convention $V = V^{+} + V^{-}$ with
$V^{+}\eqdef \max\{ V,\,0\}$ and
$V^{-}\eqdef \min\{ V,\,0\}$ (${}\leq 0$) for $V\in \RR$; nevertheless,
we will stick with our notation $V = V^{+} - V^{-}$ with
$V^{-}\eqdef \max\{ -V,\,0\}$ (${}\geq 0$).
We will not worry about this alternative any more and focus entirely
on the non\-linear equation \eqref{e:BS_M=hatV}.
Making use of eq.~\eqref{e:F(hatV)}, we arrive at
the following equivalent form of
eq.~\eqref{e:BS_M=hatV}, frequently used, cf.\
\cite[Section~2, Eq.~(1)]{Burgard-Kjaer-PDE}:
\begin{equation}
\label{eq:BS_M=hatV}
\begin{aligned}
  \frac{\partial\hat{V}}{\partial t} + \mathcal{A}_t\hat{V}
  - r\hat{V}
  = {} - (1 - R_B)\lambda_B\, \hat{V}^{-}
       + (1 - R_C)\lambda_C\, \hat{V}^{+} + s_F\, \hat{V}^{+}
\\
    \quad\mbox{ for }\, (S,t)\in (0,\infty)\times (0,T) \,.
\end{aligned}
\end{equation}
This backward parabolic equation is supplemented by
the {\it\bfseries terminal condition\/} \eqref{e:BS_M=hatV,T}.

Models with {\sl\bfseries nonlinearities\/}
are neither popular nor very frequent in Mathematical Finance.
In the present article we treat a class of semi\-linear parabolic equations
of type \eqref{eq:BS_M=hatV} with the standard linear diffusion operator
\begin{math}
  \frac{\partial}{\partial t} + \mathcal{A}_t
\end{math}
and the non\-linear reaction function that is more general that
the one on the right\--hand side of
\eqref{eq:BS_M=hatV} (only uniformly Lipschitz\--continuous).
As far as we know, this class was introduced in the work by
{\sc C.\ Burgard} and {\sc M.\ Kjaer}
\cite[Section~3]{Burgard-Kjaer-PDE} and
\cite{Burgard-Kjaer-PDE-add}.
Another class of nonlinear models is based on
a {\em\bfseries nonlinear Black\--Scholes PDE\/}
with the quasi\-linear diffusion operator
\begin{math}
    \frac{\partial\hat{V}}{\partial t}
  + \frac{1}{2}\, \sigma^2 S^2\, \frac{\partial^2 V}{\partial S^2}
  + \dots \,,
\end{math}
where the volatility
\begin{math}
  \sigma\equiv \sigma \genfrac{(}{)}{}1{\partial^2 V}{\partial S^2}
\end{math}
depends on the second partial derivative,
and with a ``typical'' linear reaction function
(sometimes including also {\em transaction costs\/}).
This class can be traced to
{\sc G.\ Barles\/} and {\sc H.~M.\ Soner\/}
\cite[Eq.\ (1.2), p.~372]{Barles-Soner}
with some additional analytic studies (on explicit solutions)
performed in
{\sc L.~A.\ Bordag\/} and {\sc Y.\ Chmakova\/}
\cite{Bordag-Chmak-2007}.
Some additional references to related numerical studies and simulations
will be added in Sections \ref{s:Numer_method} and \ref{s:Numer_parab}.

This article is organized as follows.
We begin with
a functional\--analytic reformulation of the B-S equation
\eqref{eq:BS_M=hatV}
in the next section (Section~\ref{s:Funct_anal}).
The terminal value problem \eqref{eq:BS_M=hatV}, \eqref{e:BS_M=hatV,T}
will be transformed into an initial value Cauchy problem of parabolic type.
This Cauchy problem is an initial value problem for
the non\-linear (semi\-linear) B-S equation with
a uniformly Lipschitz\--continuous (non\-linear) reaction function, as well.
In Section~\ref{s:Monotone} we construct
a {\sl\bfseries monotone iteration scheme\/}
of {\it super\-solutions\/} of this B-S equation that converge
as a monotone decreasing (i.e., non\-increasing) sequence
to the solution from above; see our main result,
Theorem~\ref{thm-Iteration}.
A closely related ramification of this monotone iteration scheme provides
an increasing sequence of {\it sub\-solutions\/} of the B-S equation
that converge to the solution from below; see
Remark~\ref{rem-Iteration}.

Numerical methods play an important role in Mathematical Finance.
In Section~\ref{s:Numer_method}
we discuss applications of two most common methods to  Mathematical Finance,
{\em finite differences/elements\/} and {\em Monte Carlo\/}.
We discuss their advantages and problems when compared to each other.
Finally, in Section~\ref{s:Numer_parab}
we derive an {\sl explicit formula\/} for the solution of
the inhomogeneous linear parabolic initial value problem for the B-S equation
that serves for computing the monotone iteration scheme in
Section~\ref{s:Monotone}.
This formula is obtained by
{\sl\bfseries variation\--of\--constants\/}
(with integrals over $\RR^1$ and $[0,T]$)
which makes it interesting for {\em Monte Carlo\/} computations.
On the other hand, the solution of the inhomogeneous linear parabolic problem
can be computed also by {\em finite differences/elements\/}.

%%%%%%%%%%%%%%%%%%%%%%%%%%%%%%%%%%%%%%%%%%%%%%%%%%%%%%%%%%%%%%%%%%%%%%%
%%%%%    Functional\--analytic reformulation    %%%%%%%%%%%%%%%%%%%%%%%
%%%%%%%%%%%%%%%%%%%%%%%%%%%%%%%%%%%%%%%%%%%%%%%%%%%%%%%%%%%%%%%%%%%%%%%

\section{Functional\--analytic reformulation of the B-S equation}
\label{s:Funct_anal}

We wish to treat the {\it terminal value problem\/}
\eqref{eq:BS_M=hatV} (or, equivalently, \eqref{e:BS_M=hatV})
above, with the terminal condition \eqref{e:BS_M=hatV,T},
by standard analytic and numerical methods for semilinear parabolic
{\it initial value problems\/}.
To this end, we rewrite problem
\eqref{eq:BS_M=hatV}, \eqref{e:BS_M=hatV,T}
as the following general initial value problem for the unknown function
$v\colon \RR^1\times  (0,T)\to \RR$,
\begin{alignat}{2}
\label{e:BS_M=v}
  \frac{\partial v}{\partial\tau} - \mathcal{A}(\tau) v
  + r\, v
& {}
  = \tilde{F}(v(x,\tau); x, \tau)
&&  \quad\mbox{ for }\, (x,\tau)\in \RR^1\times (0,T) \,;
\\
\label{e:BS_M=v,0}
    v(x,0) & {}= v_0(x)\eqdef h(\ee^x)
&&  \quad\mbox{ for }\, x\in \RR^1 \,,
\end{alignat}
where $\mathcal{A}(\tau)$ denotes
the {\it\bfseries Black\--Scholes operator\/} defined by
\begin{align}
%\label{e:BS_oper}
\label{e_tau:BS_oper}
\begin{aligned}
& (\mathcal{A}(\tau) v)(x,\tau)\eqdef
  (\mathcal{A}_{T-\tau} v)(x,\tau)
\\
& {}
  = \frac{1}{2}\, [\sigma(T-\tau)]^2 \,\frac{\partial^2 v}{\partial x^2}
  + \left( q_S(T-\tau) - \gamma_S(T-\tau)
         - \frac{1}{2}\, [\sigma(T-\tau)]^2
    \right) \frac{\partial v}{\partial x}
\end{aligned}
\\
\nonumber
    \quad\mbox{ for }\;
  v\colon \RR^1\times (0,T)\to \RR\colon
  (x,\tau)\mapsto v(x,\tau) \,,
\end{align}
and the non\-linearity
$\tilde{F}(\,\cdot\, ; x, \tau)\colon \RR\to \RR$ is given by
\begin{align}
\label{e:F(v)}
\begin{aligned}
  \tilde{F}(v; x, \tau)
& {}\eqdef
    {}- F(v; \ee^x, T-\tau) - (\lambda_B + \lambda_C)\, v
\\
& {}= (1 - R_B)\lambda_B\, v^{-}
    - (1 - R_C)\lambda_C\, v^{+} - s_F\, v^{+}
\\
\end{aligned}
\\
\nonumber
    \quad\mbox{ for }\, v\in \RR
    \,\mbox{ and }\, (x,\tau)\in \RR^1\times (0,T) \,.
\end{align}
Here, $\tau = T-t$ stands for the {\it\bfseries time to maturity\/}
and $x = \log~S$ is
the {\it\bfseries logarithmic asset (stock) price\/}; we take
$(x,\tau)\in \RR^1\times (0,T)$.
In the sequel we will never use the
{\it\bfseries real time\/} $t = T - \tau\in (0,T)$ any more,
so we prefer to use the letter $t$ in place of $\tau$ to denote
the {\it\bfseries time to maturity\/},
as it is usual in parabolic problems.
According to this new notation, in eq.~\eqref{e_tau:BS_oper}
we replace the time\--dependent coefficients
$\sigma(T-\tau)$, $q_S(T-\tau)$, and $\gamma_S(T-\tau)$ by
$\sigma(t)$, $q_S(t)$, and $\gamma_S(t)$, respectively,
and thus forget about the original terminal value problem
\eqref{eq:BS_M=hatV}, \eqref{e:BS_M=hatV,T}:
\begin{align}
\label{e:BS_oper}
\begin{aligned}
  (\mathcal{A}(t) v)(x,t)\eqdef \frac{\partial}{\partial x}\,
    \left[
    \frac{1}{2}\, [\sigma(t)]^2 \,\frac{\partial v}{\partial x}
  + \left( q_S(t) - \gamma_S(t) - \frac{1}{2}\, [\sigma(t)]^2
    \right) v(x,t)
    \right]
\end{aligned}
\\
\nonumber
    \quad\mbox{ for }\;
  v\colon \RR^1\times (0,T)\to \RR\colon
  (x,t)\mapsto v(x,t) \,.
\end{align}

Next, in order to make the initial value problem
\eqref{e:BS_M=v}, \eqref{e:BS_M=v,0}
compatible with the {\sl\bfseries monotone methods\/}
described in the article by
{\sc David H.\ Sattinger} \cite{Satting-1972},
we rewrite this problem as follows:
\begin{alignat}{2}
\label{e:BS_M=v(t)}
  \frac{\partial v}{\partial t} - \mathcal{A}(t) v
  + (r + L_{\tilde{F}})\, v
& {}
  = \tilde{F}(v(x,t); x,t) + L_{\tilde{F}}\, v
&&  \quad\mbox{ for }\, (x,t)\in \RR^1\times (0,T) \,;
\end{alignat}
with the initial condition
$v(x,0) = v_0(x)\eqdef h(\ee^x)$ for $x\in \RR^1$
in eq.~\eqref{e:BS_M=v,0},
where the constant $L_{\tilde{F}}\in \RR_+ = [0,\infty)$
is defined by
\begin{equation}
\label{e:L_tildeF}
  L_{\tilde{F}} =
  \max\{ (1 - R_B)\lambda_B ,\, (1 -  R_C)\lambda_C + s_F \} \,.
\end{equation}
According to \cite{Burgard-Kjaer-PDE},
$s_F\equiv r_F - r$, $\lambda_B\equiv r_B - r$, and
$\lambda_C\equiv r_C - r$ are some non\-negative constants and
$R_B, R_C\in [0,1]$ are the {\it\bfseries recovery rates\/}
on the derivative positions of parties $B$ and $C$, respectively.
As a consequence, the function
\begin{math}
  G(\,\cdot\, ; x,t)\colon v\mapsto G(v; x,t)\colon \RR\to \RR \,,
\end{math}
defined by
\begin{align}
\label{e:G=F-L_tildeF}
\begin{aligned}
& G(v; x,t) = \tilde{F}(v; x,t) + L_{\tilde{F}}\, v
\\
& = {}- \left[ L_{\tilde{F}} - (1 - R_B)\lambda_B \right] v^{-}
      + \left[ L_{\tilde{F}} - (1 - R_C)\lambda_C - s_F \right] v^{+}
\end{aligned}
\\
\nonumber
    \quad\mbox{ for $v\in \RR$ and }\, (x,t)\in \RR^1\times (0,T) \,,
\end{align}
is monotone increasing (i.e., non\-decreasing) on $\RR$.
Notice that both functions,
$v\mapsto {}- v^{-}$ and $v\mapsto v^{+}$,
are non\-decreasing on $\RR$.
Indeed, we have also
\begin{align*}
    \frac{\partial G}{\partial v} (v; x,t)
  = \frac{ \partial\tilde{F} }{\partial v} + L_{\tilde{F}}
  = \left\{\quad
\begin{alignedat}{2}
&   L_{\tilde{F}} - (1 - R_B)\lambda_B
&&  \quad\mbox{ if }\, v < 0 \,,
\\
&   L_{\tilde{F}} - \left( (1-R_C)\lambda_C + s_F \right)
&&  \quad\mbox{ if }\, v > 0 \,;
\end{alignedat}
\right.
\\
  \mbox{ with }\quad
    0\leq \frac{\partial G}{\partial v} (v; x,t) \leq L_{\tilde{F}}
    \quad\mbox{ for all }\, v\in \RR^1\setminus \{ 0\}
    \,\mbox{ and }\, (x,t)\in \RR\times (0,T) \,.
\end{align*}
In addition, the left\--hand side of eq.~\eqref{e:BS_M=v(t)}
clearly satisfies the weak maximum principle.

We now specify our hypotheses on the general non\-linearity
$G\colon \RR\times \RR^1\times (0,T)$
on the right\--hand side of eq.~\eqref{e:BS_M=v(t)}
treated in our present work.
We assume that $G$ satisfies the following hypotheses:

%\par\vskip 10pt
%%%%%%%%%%%%%%%%%%%%%%%%%%%%%%%%%%%%%%%%%%%%%%%%%%%%%%%%%%%%%%%%%%%%%%%
%%%%%%%%%%     General nonlinearity G(v;x,t) - HYPOTHESES (Hypos)    %%
%%%%%%%%%%%%%%%%%%%%%%%%%%%%%%%%%%%%%%%%%%%%%%%%%%%%%%%%%%%%%%%%%%%%%%%
\begin{hypos}\nopagebreak
\begingroup\rm
\begin{enumerate}
\setcounter{enumi}{0}
\renewcommand{\labelenumi}{{\bf (G\arabic{enumi})}}
\item
\makeatletter
\def\@currentlabel{{\bf G\arabic{enumi}}}\label{hy:measur}
\makeatother
For each fixed $v\in \RR$, the function\hfil\break
\begin{math}
  G(v; \,\cdot\,, \,\cdot\,)\colon (x,t)\mapsto G(v;x,t)\colon
  \RR^1\times (0,T)\to \RR
\end{math}
is Lebesgue\--measurable.
\item
\makeatletter
\def\@currentlabel{{\bf G\arabic{enumi}}}\label{hy:Lipschitz}
\makeatother
For almost every pair $(x,t)\in \RR^1\times (0,T)$, the function
\begin{math}
  G(\,\cdot\, ; x,t)\colon v\mapsto G(v; x,t)\colon \RR\to \RR
\end{math}
is uniformly Lipschitz\--continuous with a Lipschitz constant
$L_G\in \RR_+$, that is, we have
\begin{equation}
\label{e:G-Lip}
\begin{aligned}
& | G(v_1;x,t) - G(v_2;x,t) | \leq L_G\, |v_1 - v_2|
\\
&   \quad\mbox{ for all }\, v_1, v_2\in \mathbb{R}
    \,\mbox{ and for almost all }\,
    (x,t)\in \RR^1\times (0,T) \,.
\end{aligned}
\end{equation}
\item
\makeatletter
\def\@currentlabel{{\bf G\arabic{enumi}}}\label{hy:increasing}
\makeatother
For almost every pair $(x,t)\in \RR^1\times (0,T)$, the function
\begin{math}
  G(\,\cdot\, ; x,t)\colon v\mapsto G(v; x,t)\colon \RR\to \RR
\end{math}
is monotone increasing, that is,
$v_1\leq v_2$ in $\RR$ implies $G(v_1;x,t)\leq G(v_2;x,t)$.
\item
\makeatletter
\def\@currentlabel{{\bf G\arabic{enumi}}}\label{hy:exp_growth}
\makeatother
There is a constant $C_0\in \RR_+$ such that,
at almost every time $t\in (0,T)$, the function
\begin{math}
  G( 0; \,\cdot\, ,t)\colon x\mapsto G(0; x,t)\colon \RR^1\to \RR
\end{math}
satisfies the {\it exponential growth\/} restriction
\begin{equation}
\label{e:G-exp}
  | G(0;x,t) |\leq C_0\cdot \exp (|x|)
  \hspace*{10pt} \bigl( {}\leq C_0 (\ee^x + \ee^{-x}) \bigr)
    \quad\mbox{ for almost all }\, x\in \RR^1 \,.
\end{equation}
\item
\makeatletter
\def\@currentlabel{{\bf G\arabic{enumi}}}\label{hy:Hoelder_t}
\makeatother
There are constants $C_1\in \RR_+$ and $\vartheta_G\in (0,1)$ such that,
for every $v\in \RR$ and almost every $x\in \RR^1$, the function
\begin{math}
  G(v;x, \,\cdot\, )\colon t\mapsto G(v;x,t)\colon (0,T)\to \RR
\end{math}
is H{\"o}lder\--continuous with the {\it H{\"o}lder exponent\/}
$\vartheta_G$ ($0 < \vartheta_G < 1$) in the following sense,
\begin{equation}
\label{e:G-Hoelder}
\begin{aligned}
& | G(v;x,t_1) - G(v;x,t_2) |
  \leq C_1\, |v|\cdot |t_1 - t_2|^{\vartheta_G}
\\
&   \quad\mbox{ for all }\, t_1, t_2\in (0,T)
    \,\mbox{ and for almost all }\,
    (v,x)\in \RR\times \RR^1 \,.
\end{aligned}
\end{equation}
\end{enumerate}
%
%\hfill\Square
\endgroup
\end{hypos}
%%%%%%%%%%%%%%%%%%%%%%%%%%%%%%%%%%%%%%%%%%%%%%%%%%%%%%%%%%%%%%%%%%%%%%%
\par\vskip 10pt

From now on, let us consider the following generalization of
the initial value problem
\eqref{e:BS_M=v(t)}, \eqref{e:BS_M=v,0}:
\begin{alignat}{2}
\label{e:BS_v(t)}
  \frac{\partial v}{\partial t} - \mathcal{A}(t) v
  + r_G\, v
& {}
  = G(v(x,t); x,t)
&&  \quad\mbox{ for }\, (x,t)\in \RR^1\times (0,T) \,;
\end{alignat}
with the initial condition
$v(x,0) = v_0(x)\eqdef h(\ee^x)$ for $x\in \RR^1$
in eq.~\eqref{e:BS_M=v,0},
where the constant $r + L_{\tilde{F}}$ in eq.~\eqref{e:BS_M=v(t)}
has been replaced by the new constant $r_G\in \RR_+$,
owing to our monotonicity hypothesis \eqref{hy:increasing}.
Concerning hypotheses on the time\--dependent coefficients
that appear in the Black\--Scholes operator $\mathcal{A}(t)$
defined in eq.~\eqref{e:BS_oper} (recall that $\tau = t$),
we assume the following H{\"o}lder continuity:

%\par\vskip 10pt
%%%%%%%%%%%%%%%%%%%%%%%%%%%%%%%%%%%%%%%%%%%%%%%%%%%%%%%%%%%%%%%%%%%%%%%
%%%%     Hoelder continuity: Black-Scholes operator (Hypos)    %%%%%%%%
%%%%%%%%%%%%%%%%%%%%%%%%%%%%%%%%%%%%%%%%%%%%%%%%%%%%%%%%%%%%%%%%%%%%%%%
\begin{hypos}\nopagebreak
\begingroup\rm
\par\noindent
\vspace*{-10pt}
\begin{enumerate}
\setcounter{enumi}{0}
\renewcommand{\labelenumi}{{\bf (BS\arabic{enumi})}}
\item
\makeatletter
\def\@currentlabel{{\bf BS\arabic{enumi}}}\label{hy:sigma}
\makeatother
$\sigma\colon [0,T]\to (0,\infty)$
is a positive, H{\"o}lder\--continuous function satisfying
\begin{equation}
\label{e:sigma-Hoelder}
  | \sigma(t_1) - \sigma(t_2) |
  \leq C_{\sigma}\, |t_1 - t_2|^{ \vartheta_{\sigma} }
    \quad\mbox{ for all }\, t_1, t_2\in [0,T] \,,
\end{equation}
where $C_{\sigma}\in \RR_+$ and $\vartheta_{\sigma}\in (0,1)$
are some constants independent from time $t\in [0,T]$.
\item
\makeatletter
\def\@currentlabel{{\bf BS\arabic{enumi}}}\label{hy:q,gamma}
\makeatother
$q_S, \gamma_S\colon [0,T]\to \RR$
is a pair of H{\"o}lder\--continuous function satisfying
\begin{align}
\label{e:q_S:Hoelder}
  | q_S(t_1) - q_S(t_2) |
& {}
  \leq C_q\, |t_1 - t_2|^{ \vartheta_q }
    \quad\mbox{ and }\quad
\\
\label{e:gamma_S:Hoelder}
  | \gamma_S(t_1) - \gamma_S(t_2) |
& {}
  \leq C_{\gamma}\, |t_1 - t_2|^{ \vartheta_{\gamma} }
    \quad\mbox{ for all }\, t_1, t_2\in [0,T] \,,
\end{align}
where $C_q, C_{\gamma}\in \RR_+$ and
$\vartheta_q, \vartheta_{\gamma}\in (0,1)$
are some constants (independent from $t\in [0,T]$).
\end{enumerate}
%
%\hfill\Square
\endgroup
\end{hypos}
%%%%%%%%%%%%%%%%%%%%%%%%%%%%%%%%%%%%%%%%%%%%%%%%%%%%%%%%%%%%%%%%%%%%%%%
%\par\vskip 10pt

%\par\vskip 10pt
%%%%%%%%%%%%%%%%%%%%%%%%%%%%%%%%%%%%%%%%%%%%%%%%%%%%%%%%%%%%%%%%%%%%%%%
%%%%     Hoelder continuity: Black-Scholes operator (1. Remark)    %%%%
%%%%%%%%%%%%%%%%%%%%%%%%%%%%%%%%%%%%%%%%%%%%%%%%%%%%%%%%%%%%%%%%%%%%%%%
\begin{remark}\label{rem-hypos-r(t)}\nopagebreak
{\rm (H{\"o}lder exponents.)}$\;$
\begingroup\rm
In {\rm Hypotheses\/}
\eqref{hy:Hoelder_t}, \eqref{hy:sigma}, and \eqref{hy:q,gamma}
we may and will replace the {\it H{\"o}lder exponents\/}
$\vartheta_G$, $\vartheta_{\sigma}$, $\vartheta_q$,
and $\vartheta_{\gamma}$ by their minimum $\vartheta_0$,
\hfil\break
\null\hfill
\begin{math}
  \vartheta_0 =
  \min\{ \vartheta_G ,\, \vartheta_{\sigma} ,\,
         \vartheta_q ,\, \vartheta_{\gamma} \} \,,\quad
  \vartheta_0\in (0,1) \,.
\end{math}
\hfill\Square
\endgroup
\end{remark}
%%%%%%%%%%%%%%%%%%%%%%%%%%%%%%%%%%%%%%%%%%%%%%%%%%%%%%%%%%%%%%%%%%%%%%%
\par\vskip 10pt

Clearly, from {\rm Hypothesis\/} \eqref{hy:sigma} we derive
\begin{math}
  \sigma(t)\geq \sigma_0 = \min_{t\in [0,T]} \sigma(t) > 0
\end{math}
for all $t\in [0,T]$.
This fact, combined with \eqref{hy:q,gamma},
guarantees the uniform ellipticity of the Black\--Scholes operator
$\mathcal{A}(t)$ independently from $t\in [0,T]$.

%\par\vskip 10pt
%%%%%%%%%%%%%%%%%%%%%%%%%%%%%%%%%%%%%%%%%%%%%%%%%%%%%%%%%%%%%%%%%%%%%%%
%%%%     Hoelder continuity: Black-Scholes operator (2. Remark)    %%%%
%%%%%%%%%%%%%%%%%%%%%%%%%%%%%%%%%%%%%%%%%%%%%%%%%%%%%%%%%%%%%%%%%%%%%%%
\begin{remark}\label{rem-hypos-r_G}\nopagebreak
{\rm (Risk\--free interest rate.)}$\;$
\begingroup\rm
One may also suggest to replace the multiplicative constant
$r_G\in \RR$ on the left\--hand side of eq.~\eqref{e:BS_v(t)}
by the time\--dependent {\it\bfseries risk\--free interest rate\/}
$r\colon [0,T]\to \RR$ satisfying a H{\"o}lder continuity condition
analogous to those in
eqs.\ \eqref{e:q_S:Hoelder} and \eqref{e:gamma_S:Hoelder}.
However, this change would not make eq.~\eqref{e:BS_v(t)}
more general in that it could be reduced to the present form
\eqref{e:BS_v(t)} with the term $r_G\, v$ as follows:

First, define $r_G\in \RR$ by
$r_G = \max_{t\in [0,T]} r(t)$; then replace the function
$G(v;x,t)$ on the right\--hand side of eq.~\eqref{e:BS_v(t)}
by the sum
\begin{math}
  G_r(v;x,t) = G(v;x,t) + [r_G - r(t)]\, v
\end{math}
for $(v;x,t)\in \RR\times \RR^1\times (0,T)$.
Clearly, thanks to $r_G - r(t)\geq 0$ for every $t\in [0,T]$,
the function
\begin{math}
  G_r(\,\cdot\,; \,\cdot\,, \,\cdot\,)\colon
    (v;x,t)\mapsto G_r(v;x,t)\colon \RR^1\times (0,T)\to \RR
\end{math}
satisfies all {\rm Hypotheses\/} 
\eqref{hy:measur} -- \eqref{hy:Hoelder_t}
imposed on the function $G$.
We conclude that the {\it interest rate difference\/},
$r_G - r(t)$, can be included in the reaction function $G$.
We thus keep eq.~\eqref{e:BS_v(t)} in the present form with
$r_G\in \RR$ being a given constant.
\hfill\Square
\endgroup
\end{remark}
%%%%%%%%%%%%%%%%%%%%%%%%%%%%%%%%%%%%%%%%%%%%%%%%%%%%%%%%%%%%%%%%%%%%%%%
\par\vskip 10pt

Our last hypothesis in problem
\eqref{e:BS_v(t)}, \eqref{e:BS_M=v,0}
restricts the growth of the initial condition
$v(x,0) = v_0(x)\eqdef h(\ee^x)$ for $x\in \RR^1$
in eq.~\eqref{e:BS_M=v,0} as follows:

%\par\vskip 10pt
%%%%%%%%%%%%%%%%%%%%%%%%%%%%%%%%%%%%%%%%%%%%%%%%%%%%%%%%%%%%%%%%%%%%%%%
%%%%     Growth of initial values |v_0(x)| <= exp(|x|) (Hypo)    %%%%%%
%%%%%%%%%%%%%%%%%%%%%%%%%%%%%%%%%%%%%%%%%%%%%%%%%%%%%%%%%%%%%%%%%%%%%%%
\begin{hypo}\nopagebreak
\begingroup\rm
\par\noindent
\vspace*{-10pt}
\begin{enumerate}
\setcounter{enumi}{0}
\renewcommand{\labelenumi}{{\bf (v$_\mathbf{0}$)}}
\item
\makeatletter
\def\@currentlabel{{\bf v$_\mathbf{0}$}}\label{hy:v_0}
\makeatother
The function $v_0\colon \RR\to \RR$ is Lebesgue\--measurable and
there is a constant $C_h\in \RR_+$ such that,
for almost all $x\in \RR^1$, we have
\begin{equation}
\label{e:BS_v(0)}
  |v_0(x)| = |h(\ee^x)|\leq C_h\cdot \exp (|x|)
  \hspace*{10pt} \bigl( {}\leq C_h (\ee^x + \ee^{-x}) \bigr) \,.
\end{equation}
\end{enumerate}
%
%\hfill\Square
\endgroup
\end{hypo}
%%%%%%%%%%%%%%%%%%%%%%%%%%%%%%%%%%%%%%%%%%%%%%%%%%%%%%%%%%%%%%%%%%%%%%%
\par\vskip 10pt

As we have already indicated in our hypothesis \eqref{hy:exp_growth}
on the {\it exponential growth\/} restriction of $G$,
we are going to look for (strong, weak or mild) solutions
$v\colon \RR^1\times (0,T)\to \RR$ to the initial value problem
\eqref{e:BS_v(t)}, \eqref{e:BS_M=v,0}
satisfying an analogous {\it exponential growth\/} restriction of type
$v( \,\cdot\, ,t)\in H_{\CC}$ at every time $t\in (0,T)$, where
$H_{\CC} = L^2(\mathbb{R};\mathfrak{w})$
denotes the {\em complex\/} Hilbert space of
all complex\--valued Lebesgue\--measurable functions
$f\colon \RR\to \CC$ with the finite norm
\begin{equation*}
%\label{e:norm_H}
\textstyle
  \| f\|_H\eqdef
  \left( \int_{\RR} |f(x)|^2\, \mathfrak{w}(x) \,\mathrm{d}x
  \right)^{1/2} < \infty \,,
\end{equation*}
where $\mathfrak{w}(x)\eqdef \ee^{-\mu |x|}$
is a weight function with some constant $\mu\in (2,\infty)$.
This norm is induced by the inner product
\begin{equation*}
%\label{def:inn_prod}
\textstyle
  (f,g)_H\equiv
  (f,g)_{ L^2(\mathbb{R};\mathfrak{w}) } \eqdef
  \int_{\mathbb{R}} f\, \bar{g}\cdot \mathfrak{w}(x) \,\mathrm{d}x
    \quad\mbox{ for }\, f,g\in H_{\CC} \,.
\end{equation*}
As usual, the symbol $\bar{z}$ denotes
the complex conjugate of a complex number $z\in \CC$ where
$\mathbb{C} = \RR + \ii\RR$ is the complex plane.
We consider the complex Hilbert space $H_{\CC}$
only for better understanding of our applications using
{\it\bfseries holomorphic semigroups\/} in $H_{\CC}$
generated by the (unbounded) {\it Black\--Scholes operator\/}
$\mathcal{A}(t)\colon H_{\CC}\to H_{\CC}$
in eq.~\eqref{e:BS_v(t)} above.
Our solutions $v(x,t)$ to the initial value problem
\eqref{e:BS_v(t)}, \eqref{e:BS_M=v,0}
will be always real\--valued, i.e.,
$v( \,\cdot\, ,t)\in H$ at every time $t\in (0,T)$,
where $H$ denotes the closed {\em real\/} vector subspace
of all real\--valued functions $f\colon \RR\to \RR$ from $H_{\CC}$.
The domain of the differential operator $\mathcal{A}(t)$,
denoted by $\mathcal{D}(\mathcal{A}(t))$,
is a {\em complex\/} vector subspace of $H_{\CC}$
which is independent from time $t\in [0,T]$, i.e.,
\begin{math}
  \mathcal{D}(\mathcal{A}(t)) \equiv D_{\CC}\subset H_{\CC}
\end{math}
for every $t\in [0,T]$.
The vector space $D_{\CC}$ consists of all functions
$f\in H_{\CC}$ whose weak (distributional) derivatives
$f'= \frac{\mathrm{d}f}{\mathrm{d}x}$ and
$f''= \frac{\mathrm{d}^2 f}{\mathrm{d}x^2}$
belong to $H_{\CC}$, as well.
We set $D = D_{\CC}\cap H$ to denote
the closed {\em real\/} vector subspace of all real\--valued functions
$f\colon \RR\to \RR$ from $D_{\CC}$.
The vector space $D_{\CC}$ becomes a Banach space under the norm
\begin{equation*}
%\label{e:norm_D}
  \| f\|_D\eqdef \| f\|_H + \| f''\|_H
    \quad\mbox{ for }\, f\in D_{\CC} \,.
\end{equation*}
This norm is equivalent with the stronger norm
\begin{equation*}
%\label{e:norm_D}
  \vertiii{f}_D\eqdef \| f\|_H + \| f'\|_H + \| f''\|_H
    \quad\mbox{ for }\, f\in D_{\CC} \,,
\end{equation*}
by a simple interpolation inequality.
We refer to the monograph by
{\sc K.-J.\ Engel} and {\sc R.\ Nagel} \cite{Eng-Nagel}
for details concerning
{\it\bfseries holomorphic semigroups\/} and their (infinitesimal)
generators, especially to
\cite[Chapt.~II, Sect.~$4a$, pp.\ 96--109]{Eng-Nagel}.

We denote by $H_{\CC}^1$ the complex interpolation space between
$D_{\CC}$ and $H_{\CC}$ that consist of all functions
$f\colon \RR\to \RR$ from $H_{\CC}$ such that both
$f, f'\in H_{\CC}$.
$H_{\CC}^1$ is a vector space which becomes a Banach space under the norm
\begin{equation*}
%\label{e:norm_H^1}
  \| f\|_{H^1}\eqdef \| f\|_H + \| f'\|_H
    \quad\mbox{ for }\, f\in H_{\CC}^1 \,.
\end{equation*}
Hence, we have the continuous imbeddings
\begin{math}
  D_{\CC}\hookrightarrow H_{\CC}^1\hookrightarrow H_{\CC} \,.
\end{math}
Moreover, given any fixed $t\in [0,T]$,
$H_{\CC}^1$ is the domain of the {\em sesquilinear form\/}
\begin{align}
\nonumber
& \mathcal{Q}(t)\colon H_{\CC}^1\times H_{\CC}^1\to \CC\colon
  (f,g) \;\longmapsto\; \mathcal{Q}(t)(f,g)\eqdef
  {}- \int_{\mathbb{R}} [\mathcal{A}(t) f](x)\, \overline{g(x)}
                        \cdot \mathfrak{w}(x) \,\mathrm{d}x
\\
\nonumber
&
\begin{aligned}
  = {}-
& \int_{-\infty}^{+\infty}
  \left\{
          \genfrac{}{}{}1{1}{2}\, [\sigma(t)]^2\, f''(x)
  + \left( q_S(t) - \gamma_S(t)
         - \genfrac{}{}{}1{1}{2}\, [\sigma(t)]^2
    \right) f'(x)
  \right\} \overline{g(x)}\cdot \mathfrak{w}(x) \,\mathrm{d}x
\end{aligned}
\\
\label{e:sesquilin_H^1}
&
\begin{aligned}
  = \int_{-\infty}^{+\infty}
  \left\{ \genfrac{}{}{}1{1}{2}\, [\sigma(t)]^2\,
          f'(x)\, \overline{g'(x)}
  - \genfrac{}{}{}1{\mu}{2}\, [\sigma(t)]^2\,
    \mathop{\mathrm{sign}}(x)\, f'(x)\, \overline{g(x)}
  \right.
\\
  {}
  - \left.
    \left( q_S(t) - \gamma_S(t)
         - \genfrac{}{}{}1{1}{2}\, [\sigma(t)]^2
    \right) f'(x)\, \overline{g(x)}
  \right\} \cdot \mathfrak{w}(x) \,\mathrm{d}x
\end{aligned}
\\
\nonumber
&
\begin{aligned}
& {}
  = \genfrac{}{}{}1{1}{2}\, [\sigma(t)]^2 \cdot
\textstyle
    \int_{-\infty}^{+\infty}
    f'(x)\, \overline{g'(x)}\cdot \mathfrak{w}(x) \,\mathrm{d}x
\\
& {}
  + \genfrac{}{}{}1{\mu}{2}\, [\sigma(t)]^2 \cdot
\textstyle
    \left( \int_{-\infty}^0 - \int_0^{+\infty}
    \right)
    f'(x)\, \overline{g(x)}\cdot \mathfrak{w}(x) \,\mathrm{d}x
\\
& {}
\textstyle
  - \left( q_S(t) - \gamma_S(t)
         - \genfrac{}{}{}1{1}{2}\, [\sigma(t)]^2
    \right) \int_{-\infty}^{+\infty}
      f'(x)\, \overline{g(x)}\cdot \mathfrak{w}(x) \,\mathrm{d}x
\end{aligned}
\\
\nonumber
&   \quad\mbox{ defined first only for }\, f,g\in D_{\CC} \,.
\end{align}
The continuous extension of $\mathcal{Q}(t)(f,g)$
to all $f,g\in H_{\CC}^1$ is immediate, thanks to
$D_{\CC}$ being a dense vector subspace of $H_{\CC}^1$.
A few simple applications of the Cauchy\--Schwartz inequality show that
the (non\--symmetric) {\em sesquilinear form\/}
$\mathcal{Q}(t)$ on $H_{\CC}^1$ is {\em\bfseries coercive\/}.
Indeed, with a help from {\rm Hypothesis\/} \eqref{hy:sigma}
we have 
\begin{math}
  \sigma(t)\geq \sigma_0 = \min_{t\in [0,T]} \sigma(t) > 0
\end{math}
for all $t\in [0,T]$.
Consequently, if the constant $\lambda_0\in (0,\infty)$
below is chosen sufficiently large, then we get
the following more precise quantification of coercivity
at every time $t\in [0,T]$:
\begin{equation}
\label{e:Q-coercive_H^1}
  \mathcal{Q}(t)(f,f) + \lambda\, (f,f)_H
  \geq \frac{\sigma_0}{4}\, \| f\|_{H^1}^2
     + \| f\|_H^2
    \quad\mbox{ for every }\, \lambda\geq \lambda_0 \,,
\end{equation}
thanks to {\rm Hypotheses\/} \eqref{hy:sigma} and \eqref{hy:q,gamma}.
We note that the constant $\lambda_0$ depends neither on time
$t\in [0,T]$ nor on the number $\lambda\geq \lambda_0$.

Let $I\equiv I_H$ denote the identity mapping on $H_{\CC}$.
Given any real number $\lambda\geq \lambda_0$,
from ineq.~\eqref{e:Q-coercive_H^1} we infer that the linear operator
\begin{equation*}
  {}- \mathcal{A}_{\lambda}(t) \eqdef {}- \mathcal{A}(t) + \lambda I
    \colon D_{\CC}\subset H_{\CC}\to H_{\CC}
\end{equation*}
is an isomorphism of
the Banach space $D_{\CC}$ onto another Banach space $H_{\CC}$,
both, algebraically and topologically.
Now we can apply the well\--known results for abstract linear
initial value problems of parabolic type, e.g., from
 {\sc L.~C.\ Evans} \cite{Evans-98}, Chapt.~7, {\S}1.1, p.~352, or
 {\sc J.-L.\ Lions} \cite{Lions-61}, Chapt.~IV, {\S}1, p.~44, or
 \cite{Lions-71}, Chapt.~III, eq.~(1.11), p.~102,
to conclude that the inhomogeneous linear parabolic initial value problem
\begin{alignat}{2}
\label{e:BS_v(t)=f}
  \frac{\partial v}{\partial t} - \mathcal{A}(t) v
  + r_G\, v
& {}
  = f(x,t)
&&  \quad\mbox{ for }\, (x,t)\in \RR^1\times (0,T) \,;
\end{alignat}
with the initial condition
$v(x,0) = v_0(x)\eqdef h(\ee^x)$ for $x\in \RR^1$
in eq.~\eqref{e:BS_M=v,0},
possesses a unique {\it\bfseries weak solution\/} $v\colon [0,T]\to H$,
whenever the initial value $v_0\in H$ is given, such that
$v_0\colon \RR^1\to \mathbb{R}$ obeys ineq.~\eqref{e:BS_v(0)}
in {\rm Hypothesis\/} \eqref{hy:v_0}.
The weak solution $v$ is continuous as
an $H$\--valued function of time $t\in [0,T]$, that is,
$v\in C([0,T]\to H)$.
This is the ``linear part'' of the semi\-linear problem
\eqref{e:BS_v(t)}, \eqref{e:BS_M=v,0}
with a prescribed inhomogeneity $f\colon [0,T]\to H$ that is assumed to be
{\it strongly Lebesgue\--measurable\/}
% (i.e., {\it Bochner\--Lebesgue\--measurable\/})
and (essentially)
bounded on $(0,T)$, i.e.,
$f\in L^{\infty}((0,T)\to H)$.
In our case, $f\in L^{\infty}((0,T)\to H)$
follows from our stronger hypothesis below:

%\par\vskip 10pt
%%%%%%%%%%%%%%%%%%%%%%%%%%%%%%%%%%%%%%%%%%%%%%%%%%%%%%%%%%%%%%%%%%%%%%%
%%%%     NO Hoelder continuity: t --> f(t): [0,T] --> H (Hypo)    %%%%%
%%%%%%%%%%%%%%%%%%%%%%%%%%%%%%%%%%%%%%%%%%%%%%%%%%%%%%%%%%%%%%%%%%%%%%%
\begin{hypo}\nopagebreak
\begingroup\rm
\par\noindent
\vspace*{-10pt}
\begin{enumerate}
\setcounter{enumi}{0}
\renewcommand{\labelenumi}{{\bf (f\arabic{enumi})}}
\item
\makeatletter
\def\@currentlabel{{\bf f\arabic{enumi}}}\label{hy:f(t)}
\makeatother
\begin{math}
  f\colon [0,T]\to H\subset H_{\CC} = L^2(\mathbb{R};\mathfrak{w})
\end{math}
is a continuous function, i.e., $f\in C([0,T]\to H)$, where
\begin{math}
  f\colon \RR^1\times [0,T]\to \RR\colon (x,t)\mapsto f(x,t)
\end{math}
satisfies
$f(t)\equiv f( \,\cdot\, ,t)\in H$ for every $t\in [0,T]$.
\end{enumerate}
%
%\hfill\Square
\endgroup
\end{hypo}
%%%%%%%%%%%%%%%%%%%%%%%%%%%%%%%%%%%%%%%%%%%%%%%%%%%%%%%%%%%%%%%%%%%%%%%
\par\vskip 10pt

Although in the following sections we work only with weak solutions
$v\colon [0,T]\to H$, $v(0) = v_0\in H$,
to the ``linear part'' \eqref{e:BS_v(t)=f}
of the semi\-linear equation \eqref{e:BS_v(t)},
we would like to remark that 
the unique {\it\bfseries weak solution\/}
$v\colon [0,T]\to H$ to the linear initial value problem
\eqref{e:BS_v(t)=f}, \eqref{e:BS_M=v,0}
becomes a (unique) {\it\bfseries classical solution\/} if 
$f\colon [0,T]\to H$ satisfies the following stronger,
H{\"o}lder\--continuity hypothesis
with the {\it H{\"o}lder exponent\/}
$\vartheta_f\in (0,1)$:

%\newpage
%\vfill\eject
%\par\vskip 10pt
%%%%%%%%%%%%%%%%%%%%%%%%%%%%%%%%%%%%%%%%%%%%%%%%%%%%%%%%%%%%%%%%%%%%%%%
%%%%     Hoelder continuity: t --> f(t): [0,T] --> H (Hypo)    %%%%%%%%
%%%%%%%%%%%%%%%%%%%%%%%%%%%%%%%%%%%%%%%%%%%%%%%%%%%%%%%%%%%%%%%%%%%%%%%
\begin{hypo}\nopagebreak
\begingroup\rm
\par\noindent
\vspace*{-10pt}
\begin{enumerate}
\setcounter{enumi}{0}
\renewcommand{\labelenumi}{{\bf (f\arabic{enumi}')}}
\item
\makeatletter
\def\@currentlabel{{\bf f\arabic{enumi}}}\label{hy:f-Hoelder}
\makeatother
\begin{math}
  f\colon [0,T]\to H\subset H_{\CC} = L^2(\mathbb{R};\mathfrak{w})
\end{math}
is a $\vartheta_f$\--H{\"o}lder\--continuous function, i.e.,
there are constants $\vartheta_f\in (0,1)$ and $C_f\in \RR_+$ such that
\begin{equation}
\label{e:f-Hoelder}
  \| f( \,\cdot\, ,t_1) - f( \,\cdot\, ,t_2)\|_H
  \leq C_f\cdot |t_1 - t_2|^{\vartheta_f}
    \quad\mbox{ for all }\, t_1, t_2\in [0,T] \,.
\end{equation}
\end{enumerate}
%
%\hfill\Square
\endgroup
\end{hypo}
%%%%%%%%%%%%%%%%%%%%%%%%%%%%%%%%%%%%%%%%%%%%%%%%%%%%%%%%%%%%%%%%%%%%%%%
\par\vskip 10pt

This is the case if the inhomogeneity
$f\colon \RR^1\times [0,T]\to \RR$ satisfies
$f( \,\cdot\, ,t)\in H$ for each $t\in [0,T]$
and there are some constants
$\tilde{C}_f\in \RR_+$ and $\kappa\in \RR$, with $1\leq \kappa < \mu/2$,
such that
\begin{equation}
\label{e:f(x,t)-Hoelder}
\begin{aligned}
& |f(x,t_1) - f(x,t_2)|
  \leq \tilde{C}_f\, \ee^{\kappa |x|}\cdot |t_1 - t_2|^{\vartheta_f}
\\
&   \quad\mbox{ for a.e.\ (almost every) $x\in \RR^1$ and for all }\,
    t_1, t_2\in [0,T] \,.
\end{aligned}
\end{equation}
Recall that $\mu$ ($\mu > 2$) is the constant in the weight function
$\mathfrak{w}(x)\eqdef \ee^{-\mu |x|}$ in the Hilbert space
$H_{\CC} = L^2(\mathbb{R};\mathfrak{w})$.
We remark that for the inhomogeneous linear equation
\eqref{e:BS_v(t)=f},
the nonlinearity $G$ on the right\--hand side of eq.~\eqref{e:BS_v(t)}
becomes $G(v;x,t)\equiv f(x,t)$; hence, we have
$\vartheta_G = \vartheta_f$ in {\rm Hypothesis\/} \eqref{hy:Hoelder_t}.
Let us recall that, by Remark~\ref{rem-hypos-r(t)},
we have replaced the {\it H{\"o}lder exponents\/}
$\vartheta_G$, $\vartheta_{\sigma}$, $\vartheta_q$,
and $\vartheta_{\gamma}$ by their minimum $\vartheta_0$;
hence, we may include also the value of $\vartheta_f$ in that minimum:
\begin{equation}
\label{e:min-Hoelder}
  \vartheta_0 =
  \min\{ \vartheta_G ,\, \vartheta_{\sigma} ,\,
         \vartheta_q ,\, \vartheta_{\gamma} ,\, \vartheta_f \} \,,\quad
  \vartheta_0\in (0,1) \,. 
\end{equation}

Indeed, according to the existence, uniqueness, and regularity results for
problem \eqref{e:BS_v(t)=f}, \eqref{e:BS_M=v,0}
in {\sc A.\ Pazy} \cite[Chapt.~5, {\S}5.7]{Pazy}, Theorem 7.1 on p.~168,
if \eqref{e:f-Hoelder} holds, then
the unique {\it\bfseries weak solution\/}
$v\colon [0,T]\to H$ to the linear initial value problem
\eqref{e:BS_v(t)=f}, \eqref{e:BS_M=v,0}
described above happens to be a unique {\it\bfseries classical solution\/}
which, among other properties, is continuous as a function
$v\colon [0,T]\to H$, i.e., $v\in C([0,T]\to H)$,
continuously differentiable on the time interval $(0,T]$,
$v(t)\equiv v( \,\cdot\, ,t)\in D$ for every $t\in (0,T]$, i.e.,
$\mathcal{A}(t) v(t)\in H$ for $t\in (0,T]$,
and $v$ satisfies the abstract differential equation
\begin{alignat}{2}
\label{e:abstr:BS_v(t)=f}
  \frac{\partial v}{\partial t} - \mathcal{A}(t) v
  + r_G\, v
& {}
  = f(t)
&&  \quad\mbox{ in $H$ for }\, t\in (0,T) \,;
\\
\label{e:abstr:BS_v_0}
    \quad\mbox{ with }\quad
    v(0) & {}= v_0\in H
&&  \quad\mbox{ (the initial condition).}
\end{alignat}
%

%%%%%%%%%%%%%%%%%%%%%%%%%%%%%%%%%%%%%%%%%%%%%%%%%%%%%%%%%%%%%%%%%%%%%%%
%%%%%    Monotone methods for parabolic problems    %%%%%%%%%%%%%%%%%%%
%%%%%%%%%%%%%%%%%%%%%%%%%%%%%%%%%%%%%%%%%%%%%%%%%%%%%%%%%%%%%%%%%%%%%%%

\section{Monotone methods for the non\-linear B-S equation}
\label{s:Monotone}

We make use of the inhomogeneous linear problem
\eqref{e:BS_v(t)=f}, \eqref{e:BS_M=v,0}
in order to describe an iterative scheme for approximating
the unique weak solution
$v\colon [0,T]\to H$, $v(0) = v_0\in H$,
to the semi\-linear problem
\eqref{e:BS_v(t)}, \eqref{e:BS_M=v,0}.

%%%%%%%%%%%%%%%%%%%%%%%%%%%%%%%%%%%%%%%%%%%%%%%%%%%%%%%%%%%%%%%%%%%%%%%
%%%%%    Preliminary comparison results    %%%%%%%%%%%%%%%%%%%%%%%%%%%%
%%%%%%%%%%%%%%%%%%%%%%%%%%%%%%%%%%%%%%%%%%%%%%%%%%%%%%%%%%%%%%%%%%%%%%%

\subsection{Preliminary comparison results for parabolic problems}
\label{ss:Comparison}

First, the so\--called {\it\bfseries weak maximum principle\/}
for a {\it classical solution\/} $v$
of the inhomogeneous linear problem
\eqref{e:BS_v(t)=f}, \eqref{e:BS_M=v,0}
is established e.g.\ in
{\sc A.\ Friedman\/}
\cite[Chapt.~2, Sect.~4, Theorem~9, p.~43]{Friedman-64}.
A standard approximation procedure of a {\it weak solution\/} $v$
by a sequence of classical solutions yields
the corresponding {\it weak maximum principle\/} also for
the {\it weak solution\/} $v$.
More precisely, if the inequalities
$v(x,0) = v_0(x)\geq 0$ and $f(x,t)\geq 0$
are valid for almost all $(x,t)\in \RR^1\times (0,T)$,
then also $v(x,t)\geq 0$ holds for almost all
$(x,t)\in \RR^1\times (0,T)$.

Second, let $v\colon [0,T]\to H$ be a {\it classical solution\/}
of the inhomogeneous linear problem
\begin{alignat}{2}
\label{e:abstr:BS_v(t)=g}
  \frac{\partial v}{\partial t} - \mathcal{A}(t) v
  + r_G\, v
& {}
  = g(t)
&&  \quad\mbox{ in $H$ for }\, t\in (0,T) \,;
\\
\label{e:abstr:BS_v_g}
    \quad\mbox{ with }\quad
    v(0) & {}= v_g\in H
&&  \quad\mbox{ (the initial condition),}
\end{alignat}
where $g\colon [0,T]\to H$ is a function continuous and bounded in
$(0,T)$, i.e.,
$g\in C((0,T)\to H)$ with
\begin{math}
  \| g\|_{ L^{\infty}((0,T)\to H) }\eqdef
  \sup_{t\in (0,T)} \| g(t)\|_H < \infty .
\end{math}
We say that $v$ is a {\it\bfseries super\-solution\/}
to the inhomogeneous linear problem
\eqref{e:BS_v(t)=f}, \eqref{e:BS_M=v,0},
if the inequalities $v_g(x)\geq v_0(x)$ and $g(x,t)\geq f(x,t)$
hold for almost all $(x,t)\in \RR^1\times (0,T)$.
Analogously, a {\it\bfseries sub\-solution\/} $v$
to the inhomogeneous linear problem
\eqref{e:BS_v(t)=f}, \eqref{e:BS_M=v,0}
is a {\it classical solution\/} $v$ of problem
\eqref{e:abstr:BS_v(t)=g}, \eqref{e:abstr:BS_v_g}
for which the inequalities $v_g(x)\leq v_0(x)$ and $g(x,t)\leq f(x,t)$
hold for almost all $(x,t)\in \RR^1\times (0,T)$.
More generally, if
$v\colon (x,t)\in \RR^1\times (0,T)\to \RR$
is a {\it weak solution\/} to the inhomogeneous linear problem
\eqref{e:BS_v(t)=f}, \eqref{e:BS_M=v,0},
then the notions of
{\it\bfseries super-\/} and {\it\bfseries sub\-solution\/}
to the inhomogeneous linear problem
\eqref{e:BS_v(t)=f}, \eqref{e:BS_M=v,0}
are defined by means of an approximation procedure
by a sequence of classical solutions again.
A rigorous functional\--analytic way to define
{\it\bfseries super-\/} and {\it\bfseries sub\-solution\/}
$v$ to problem
\eqref{e:BS_v(t)=f}, \eqref{e:BS_M=v,0}
is to require that
$v\colon (x,t)\in \RR^1\times (0,T)\to \RR$
have all regularity properties of a {\it weak solution\/} to problem
\eqref{e:BS_v(t)=f}, \eqref{e:BS_M=v,0}
stated in 
 {\sc L.~C.\ Evans} \cite{Evans-98}, Chapt.~7, {\S}1.1, p.~352, or
 {\sc J.-L.\ Lions} \cite{Lions-61}, Chapt.~IV, {\S}1, p.~44, or
 \cite{Lions-71}, Chapt.~III, eq.~(1.11), p.~102,
and, in addition to these regularity properties,
the following inequality is valid in the {\em sense of distributions\/} on
$\RR^1\times (0,T)$:
\begin{alignat}{2}
\label{ineq:abstr:BS_v(t)=f}
  \frac{\partial v}{\partial t} - \mathcal{A}(t) v
  + r_G\, v
& {}
  \geq f(t) \quad ({}\leq f(t))
&&  \quad\mbox{ in $H$ for }\, t\in (0,T) \,;
\\
\label{ineq:abstr:BS_v_f}
    \quad\mbox{ with }\quad
    v(0) & {}\geq v_0\in H \quad
           ({}\leq v_0\in H)
&&  \quad\mbox{ (the initial condition).}
\end{alignat}
Here, the inequalities with ``$\,\geq\,$''
(``$\,\leq\,$'', respectively)
specify a (weak) {\it\bfseries super\-solution\/}
(a (weak) {\it\bfseries sub\-solution\/}).
The reader is referred to
{\sc A.\ Friedman\/}
\cite[Chapt.~3, Sect.~3, Theorem~8, p.~51]{Friedman-63}
for details about positive (or non\-negative) distributions.
Clearly, any function
$v\colon (x,t)\in \RR^1\times (0,T)\to \RR$
which is simultaneously
a (weak) {\it\bfseries super\-solution\/} and
a (weak) {\it\bfseries sub\-solution\/})
of the inhomogeneous linear problem
\eqref{e:BS_v(t)=f}, \eqref{e:BS_M=v,0}
is a (weak) solution to this problem.
Combining these definitions of
{\it\bfseries super-\/} and {\it\bfseries sub\-solution\/},
denoted by
$\overline{v}, \underline{v}\colon \RR^1\times (0,T)\to \RR$,
respectively, having the initial values satisfying
$\underline{v}(0)\leq v_0\leq \overline{v}(0)$
a.e.\ in $\RR^1$, with the weak maximum principle for
the difference $w = \overline{v} - \underline{v}$,
we obtain the following auxiliary {\em\bfseries weak comparison\/} result.

%%%%%%%%%%%%%%%%%%%%%%%%%%%%%%%%%%%%%%%%%%%%%%%%%%%%%%%%%%%%%%%%%%%%%%%
%%%%%    Sub-solution <= Sub-solution (Lemma)    %%%%%%%%%%%%%%%%%%%%%%
%%%%%%%%%%%%%%%%%%%%%%%%%%%%%%%%%%%%%%%%%%%%%%%%%%%%%%%%%%%%%%%%%%%%%%%
\begin{lemma}[Weak comparison.]
\label{lem-Sub<=Supersol}
$\;$
Assume that\/
$\overline{v}, \underline{v}\colon \RR^1\times (0,T)\to \RR$,
respectively, is a pair of (weak)\/
{\it\bfseries super-\/} and {\it\bfseries sub\-solutions\/}
of problem \eqref{e:BS_v(t)=f}, \eqref{e:BS_M=v,0} satisfying
$\underline{v}(0)\leq \overline{v}(0)$
a.e.\ in $\RR^1$.
Then, at every time\/ $t\in [0,T)$, we have\/
$\underline{v}(x,t)\leq \overline{v}(x,t)$
for a.e.\ $x\in \RR^1$.
\end{lemma}
%%%%%%%%%%%%%%%%%%%%%%%%%%%%%%%%%%%%%%%%%%%%%%%%%%%%%%%%%%%%%%%%%%%%%%%
\par\vskip 10pt

Observe that we have left the initial value $v_0\in H$
out of this lemma since we use it usually with either
$\underline{v}(0) = v_0$ or $\overline{v}(0) = v_0$
as the initial condition attached to the differential equation
\begin{alignat}{2}
\label{e:BS_v(t)=f_sub}
  \frac{ \partial\underline{v} }{\partial t}
  - \mathcal{A}(t) \underline{v} + r_G\, \underline{v}
& {}
  = \underline{f}(x,t)
&&  \quad\mbox{ for }\, (x,t)\in \RR^1\times (0,T) \,,
\\
\intertext{or}
%\\
\label{e:BS_v(t)=f_sup}
  \frac{ \partial\overline{v} }{\partial t}
  - \mathcal{A}(t) \overline{v} + r_G\, \overline{v}
& {}
  = \overline{f}(x,t)
&&  \quad\mbox{ for }\, (x,t)\in \RR^1\times (0,T) \,,
\end{alignat}
respectively, where
\begin{math}
  \underline{f}( \,\cdot\, ,t)\leq \overline{f}( \,\cdot\, ,t)
\end{math}
a.e.\ in $\RR^1$, at every time $t\in (0,T)$.

%\par\vskip 10pt
%%%%%%%%%%%%%%%%%
%\proof
{\it Proof of\/} {\bf Lemma~\ref{lem-Sub<=Supersol}.}$\;$
We subtract equation
\eqref{e:BS_v(t)=f_sub} from \eqref{e:BS_v(t)=f_sup},
thus obtaining an analogous equation for the difference
$w = \overline{v} - \underline{v}$ with the right\--hand side equal to
$g(x,t) = \overline{f}(x,t) - \underline{f}(x,t)\geq 0$
for a.e.\ $(x,t)\in \RR^1\times (0,T)$.
Then the desired result,
$w( \,\cdot\, ,t)\geq 0$ a.e.\ in $\RR^1$,
at every time\/ $t\in [0,T)$, follows from
{\sc A.\ Friedman\/}
\cite[Chapt.~2, Sect.~4, Theorem~9, p.~43]{Friedman-64}, cf.\ also
\cite[Chapt.~2, Sect.~6, Theorem~16, p.~52]{Friedman-64}.
%\null\hfill\qed
\qed
%%%%%%%%%%%%%%%%%
\par\vskip 10pt

We now give simple {\bf examples} of super- and sub\-solutions
of problem \eqref{e:BS_v(t)}, \eqref{e:BS_M=v,0}.

%\par\vskip 10pt
%%%%%%%%%%%%%%%%%%%%%%%%%%%%%%%%%%%%%%%%%%%%%%%%%%%%%%%%%%%%%%%%%%%%%%%
%%%%%    Sub-solution <= Sub-solution (Example)    %%%%%%%%%%%%%%%%%%%%
%%%%%%%%%%%%%%%%%%%%%%%%%%%%%%%%%%%%%%%%%%%%%%%%%%%%%%%%%%%%%%%%%%%%%%%
\begin{example}[{\sl Super-\/} and {\sl sub\-solutions\/}.]
\label{exam-Sub<=Supersol}\nopagebreak
\begingroup\rm
$\;$
Let us define the function
\begin{equation}
\label{e:Supersol}
  V(x,t) = K\, \ee^{\lambda t}
           \left( \ee^{\kappa x} + \ee^{-\kappa x} \right)
         = 2K\, \ee^{\lambda t}\cdot \cosh(\kappa x)
    \quad\mbox{ for }\, (x,t)\in \RR^1\times [0,T] \,,
\end{equation}
where $\kappa\in \RR$ is a constant satisfying $1\leq \kappa < \mu/2$,
and $K, \lambda\in \RR$ with $K\geq 1$ and $\lambda\geq 0$
are some other constants (large enough)
to be determined below:

The left\--hand side of eq.~\eqref{e:BS_v(t)} with $v=V$ becomes
\begin{equation}
\label{e_left:BS_v=V(t)}
\begin{aligned}
& \mathrm{l.h.s.}(x,t) =
    \frac{\partial V}{\partial t} - \mathcal{A}(t) V + r_G\, V(x,t)
\\
& {}
  = \lambda\, V(x,t) - \frac{1}{2}\, \kappa^2 [\sigma(t)]^2\, V(x,t)
\\
& {}
  - \kappa
    \left[ q_S(t) - \gamma_S(t) - \frac{1}{2}\, [\sigma(t)]^2\right]
  \cdot
  \frac{ \ee^{\kappa x} - \ee^{-\kappa x} }%
       { \ee^{\kappa x} + \ee^{-\kappa x} }\, V(x,t)
\\
& \geq \left[ \lambda - \frac{1}{2}\, \kappa^2 [\sigma(t)]^2\right] V(x,t)
  - \kappa\cdot
    \left\vert q_S(t) - \gamma_S(t) - \frac{1}{2}\, [\sigma(t)]^2
    \right\vert\, V(x,t)
%    \pmb{\bigg\vert} q_S(t) - \gamma_S(t) - \frac{1}{2}\, [\sigma(t)]^2
%    \pmb{\bigg\vert}\, V(x,t)
\\
& \geq \left\{ \lambda - \genfrac{}{}{}1{1}{2}\,
               \kappa^2 \|\sigma\|_{ L^{\infty}(0,T) }^2
  - \kappa
    \left[ \| q_S - \gamma_S\|_{ L^{\infty}(0,T) }
         + \genfrac{}{}{}1{1}{2}\, \|\sigma\|_{ L^{\infty}(0,T) }^2
    \right] \right\} V(x,t)
\\
&   \quad\mbox{ for }\, (x,t)\in \RR^1\times [0,T] \,.
\end{aligned}
\end{equation}
As usual, $\|\sigma\|_{ L^{\infty}(0,T) }$
stands for the supremum norm of a continuous function
$\sigma\colon [0,T]\to \RR$.

On the other hand, the right\--hand side of eq.~\eqref{e:BS_v(t)}
with $v=V$ becomes
\begin{align}
\label{e_right:BS_v=V(t)}
&
\begin{aligned}
& \mathrm{r.h.s.}(x,t) = G(V(x,t); x,t) = G(0; x,t)
  + \left[ G(V(x,t); x,t) - G(0; x,t)\right]
\\
& {}
  \leq C_0\, \exp (|x|) + L_G\, V(x,t)\leq C_0'\, V(x,t)
    \quad\mbox{ for }\, (x,t)\in \RR^1\times [0,T] \,,
\end{aligned}
\\
%\label{e_right:(C_0/K)+L_G}
\nonumber
&   \quad\mbox{ where }\qquad
  C_0'\eqdef ({C_0}/{K}) + L_G \hspace*{10pt} ({} > 0) \,,
\end{align}
and we have taken advantage of inequalities
\eqref{e:G-Lip} and \eqref{e:G-exp}
in {\rm Hypotheses\/} \eqref{hy:Lipschitz} and \eqref{hy:exp_growth},
respectively.
Subtracting eq.~\eqref{e_right:BS_v=V(t)} from \eqref{e_left:BS_v=V(t)}
we arrive at
\begin{align}
\nonumber
& \mathrm{l.h.s.}(x,t) - \mathrm{r.h.s.}(x,t) =
    \frac{\partial V}{\partial t} - \mathcal{A}(t) V + r_G\, V(x,t)
  - G(V(x,t); x,t)
\\
\label{e:left-right:V(t)}
&
\begin{aligned}
  \geq \Bigl\{ \lambda
& {}
  - \genfrac{}{}{}1{1}{2}\,
    \kappa^2 \|\sigma\|_{ L^{\infty}(0,T) }^2
\\
& {}
  - \kappa
    \left[ \| q_S - \gamma_S\|_{ L^{\infty}(0,T) }
         + \genfrac{}{}{}1{1}{2}\, \|\sigma\|_{ L^{\infty}(0,T) }^2
    \right]
  - C_0'\Bigr\} V(x,t)\geq 0
\end{aligned}
\\
\nonumber
&       \quad\mbox{ for }\, (x,t)\in \RR^1\times [0,T] \,,
\end{align}
provided $\lambda\in \RR$ satisfies
\begin{equation}
\label{e:left-right:lambda}
\begin{aligned}
  \lambda\geq \Lambda\eqdef
& {}
  \genfrac{}{}{}1{1}{2}\,
    \kappa^2 \|\sigma\|_{ L^{\infty}(0,T) }^2
\\
& {}
  + \kappa
    \left[ \| q_S - \gamma_S\|_{ L^{\infty}(0,T) }
         + \genfrac{}{}{}1{1}{2}\, \|\sigma\|_{ L^{\infty}(0,T) }^2
    \right]
  + C_0'\quad ({} > 0) \,.
\end{aligned}
\end{equation}
Recalling our {\rm Hypothesis\/} \eqref{hy:v_0} with ineq.~\eqref{e:BS_v(0)}
on the growth of the initial condition,
we take the constant $K\in \RR$ such that
$K\geq \max\{ 1,\, C_h\}$ which guarantees also
\begin{equation}
\label{e:BS_V(0)}
\begin{aligned}
& |v_0(x)| = |h(\ee^x)|\leq C_h\cdot \exp (|x|)
  \leq C_h\cdot \exp (\kappa |x|)
\\
& {}
  \leq
  V(x,0) = K\, \left( \ee^{\kappa x} + \ee^{-\kappa x} \right)
         = 2K\cdot \cosh(\kappa x)
    \quad\mbox{ for }\, (x,t)\in \RR^1\times [0,T] \,.
\end{aligned}
\end{equation}
It follows that the function
$V\colon \RR^1\times [0,T]\to \RR$ defined in eq.~\eqref{e:Supersol}
is a {\em super\-solution\/}
of problem \eqref{e:BS_v(t)}, \eqref{e:BS_M=v,0}.

Analogous arguments show that the function
${}- V\colon \RR^1\times [0,T]\to \RR$
is a {\em sub\-solution\/}
of problem \eqref{e:BS_v(t)}, \eqref{e:BS_M=v,0}.
Notice that in this case, ineq.~\eqref{e_right:BS_v=V(t)}
has to be replaced by
\begin{equation}
\label{e_right:BS_v=-V(t)}
\begin{aligned}
& G(-V(x,t); x,t) = G(0; x,t)
  + \left[ G(-V(x,t); x,t) - G(0; x,t)\right]
\\
& {}
  \geq {}- C_0'\, V(x,t)
    \quad\mbox{ for }\, (x,t)\in \RR^1\times [0,T] \,.
\end{aligned}
\end{equation}
\null\hfill\Square
%\hfill\Square
\endgroup
\end{example}
%%%%%%%%%%%%%%%%%%%%%%%%%%%%%%%%%%%%%%%%%%%%%%%%%%%%%%%%%%%%%%%%%%%%%%%
%\par\vskip 10pt

%%%%%%%%%%%%%%%%%%%%%%%%%%%%%%%%%%%%%%%%%%%%%%%%%%%%%%%%%%%%%%%%%%%%%%%
%%%%%    Preliminary comparison results    %%%%%%%%%%%%%%%%%%%%%%%%%%%%
%%%%%%%%%%%%%%%%%%%%%%%%%%%%%%%%%%%%%%%%%%%%%%%%%%%%%%%%%%%%%%%%%%%%%%%

\subsection{Construction of monotone iterations}
\label{ss:Iteration}

Recalling Lemma~\ref{lem-Sub<=Supersol},
our {\rm Hypothesis\/} \eqref{hy:v_0} with ineq.~\eqref{e:BS_v(0)},
and applying Example~\ref{exam-Sub<=Supersol} with $\kappa = 1$
$(< \mu/2)$,
we are now ready to construct a {\sl\bfseries monotone iteration scheme\/}
for calculating a (weak) solution
$v\colon \RR^1\times (0,T)\to \RR$ to the initial value problem
\eqref{e:BS_v(t)}, \eqref{e:BS_M=v,0}.
We start by setting $\kappa = 1$ and fixing the constants
$K, \lambda\in \RR$ with $K\geq 1$ and
$\lambda\geq \Lambda$ $(\, > 0)$ large enough,
such that both inequalities,
\eqref{e:left-right:lambda} and \eqref{e:BS_V(0)}, are valid.
It follows from eq.~\eqref{e:Supersol}
and ineq.~\eqref{e:left-right:V(t)} that the function
$u_0\colon \RR^1\times (0,T)\to \RR$ defined by
\begin{equation}
\label{e:u_0:Supersol}
  u_0(x,t) = K\, \ee^{\lambda t} \left( \ee^x + \ee^{- x} \right)
           = 2K\, \ee^{\lambda t}\cdot \cosh(x)
    \quad\mbox{ for }\, (x,t)\in \RR^1\times [0,T] \,,
\end{equation}
is a {\em supersolution\/} of problem
\eqref{e:BS_v(t)}, \eqref{e:BS_M=v,0}.
We remark that ${}- u_0$ happens to be a {\em subsolution\/}
of this problem, by Example~\ref{exam-Sub<=Supersol} with $\kappa = 1$,
as well.
The first iterate,
$u_1\colon \RR^1\times [0,T]\to \RR$,
is constructed as the (weak) solution $u_1$ to the following analogue of
the inhomogeneous linear initial value problem
\eqref{e:BS_v(t)=f}, \eqref{e:BS_M=v,0}:
\begin{alignat}{2}
\label{e:BS_u_1(t)=f}
  \frac{\partial u_1}{\partial t} - \mathcal{A}(t) u_1
  + r_G\, u_1
& {}
  = G(u_0(x,t); x,t)
\\
\nonumber
& \left( {}\leq C_0'\, u_0(x,t) \right)
&&  \quad\mbox{ for }\, (x,t)\in \RR^1\times [0,T] \,;
\\
\label{e:BS_M=u_1,0}
    u_1(x,0) & {}= v_0(x)\eqdef h(\ee^x)
&&  \quad\mbox{ for }\, x\in \RR^1 \,.
\end{alignat}
Since $u_0$ is a (weak) supersolution of problem
\eqref{e:BS_v(t)}, \eqref{e:BS_M=v,0},
by Example~\ref{exam-Sub<=Supersol} with $\kappa = 1$,
we may apply Lemma~\ref{lem-Sub<=Supersol}
to conclude that $u_1(x,t)\leq u_0(x,t)$ holds for a.e.\
$(x,t)\in \RR^1\times (0,T)$.
In addition, making use of eq.~\eqref{e_right:BS_v=-V(t)},
we get also
\begin{equation}
\label{e:u_1<=u_0}
\begin{aligned}
  {}- u_0(x,t)\leq u_1(x,t)\leq u_0(x,t)
  \quad \left( {} = K\, \ee^{\lambda t}
          \left( \ee^{\kappa x} + \ee^{-\kappa x} \right)
        \right)
\\
    \quad\mbox{ for }\, (x,t)\in \RR^1\times [0,T] \,.
\end{aligned}
\end{equation}

Our next step is the following {\sl induction hypothesis\/}.
Let us assume that, for some integer $m\geq 1$,
in addition to $u_0$ and $u_1$ above,
we have already constructed the first $(m+1)$ functions
\begin{math}
  u_0, u_1, u_2,\dots, u_m\colon \RR^1\times [0,T]\to \RR
\end{math}
with the following properties:
\begin{itemize}
\item[{\rm (a)}]
Every function
$u_j\colon \RR^1\times [0,T]\to \RR$; $j = 0,1,2,\dots, m$,
is Lebesgue\--measurable and continuous in time as a function
$u_j\colon [0,T]\to H$, i.e., $u_j\in C([0,T]\to H)$.
\item[{\rm (b)}]
The inequalities
\begin{equation}
\label{e:u_j<=u_0}
\begin{aligned}
  {}- u_0(x,t)\leq u_j(x,t)\leq u_{j-1}(x,t)\leq u_0(x,t)
    \quad\mbox{ for }\, (x,t)\in \RR^1\times [0,T] \,,
\end{aligned}
\end{equation}
are valid for every $j = 1,2,3,\dots, m$.
\item[{\rm (c)}]
For each $j = 1,2,3,\dots, m$, the function
$u_j\colon \RR^1\times [0,T]\to \RR$
is the (weak) solution to the following analogue of
the inhomogeneous linear initial value problem
\eqref{e:BS_v(t)=f}, \eqref{e:BS_M=v,0}; cf.\ problem
\eqref{e:BS_u_1(t)=f}, \eqref{e:BS_M=u_1,0} above:
\begin{alignat}{2}
\label{e:BS_u_j(t)=f}
  \frac{\partial u_j}{\partial t} - \mathcal{A}(t) u_j
  + r_G\, u_j
& {}
  = G(u_{j-1}(x,t); x,t)
\\
\nonumber
& \left( {}\leq C_0'\, u_0(x,t) \right)
&&  \quad\mbox{ for }\, (x,t)\in \RR^1\times [0,T] \,;
\\
\label{e:BS_M=u_j,0}
    u_j(x,0) & {}= v_0(x)\eqdef h(\ee^x)
&&  \quad\mbox{ for }\, x\in \RR^1 \,.
\end{alignat}
\end{itemize}

In our last step
({\sl induction on the index $m\geq 1$\/})
we construct the $(m+1)$\--st iterate,
$u_{m+1}\colon \RR^1\times [0,T]\to \RR$,
to be the (weak) solution $u_{m+1}$ to the following problem; cf.\ problem
\eqref{e:BS_u_1(t)=f}, \eqref{e:BS_M=u_1,0}:
\begin{alignat}{2}
\label{e:BS_u_m+1(t)=f}
  \frac{\partial u_{m+1}}{\partial t}
& {}
  - \mathcal{A}(t) u_{m+1} + r_G\, u_{m+1}
&&{}
  = G(u_m(x,t); x,t)
\\
\nonumber
& \left( {}\leq C_0'\, u_0(x,t) \right)
&&  \mbox{ for }\, (x,t)\in \RR^1\times [0,T] \,;
\\
\label{e:BS_M=u_m+1,0}
    u_{m+1}(x,0) & {}= v_0(x)\eqdef h(\ee^x)
&&  \mbox{ for }\, x\in \RR^1 \,.
\end{alignat}
By arguments analogous to those used in the construction of $u_1$
from $u_0$ above, we conclude that $u_{m+1}$ exists and satisfies
\begin{equation}
\label{e:u_m+1<=u_m}
\begin{aligned}
  {}- u_0(x,t)\leq u_{m+1}(x,t)\leq u_m(x,t)\leq u_0(x,t)
    \quad\mbox{ for }\, (x,t)\in \RR^1\times [0,T] \,.
\end{aligned}
\end{equation}
Here, we have used our monotonicity hypothesis \eqref{hy:increasing}
to conclude that
$u_m\leq u_{m-1}$ a.e.\ in $\RR^1\times [0,T]$ entails
$G(u_m(x,t); x,t)\leq G(u_{m-1}(x,t); x,t)$ for a.e.\
$(x,t)\in \RR^1\times [0,T]$.
Finally, we get $u_{m+1}\in C([0,T]\to H)$.
This concludes the construction of the desired iterates.

%\par\vskip 10pt
%%%%%%%%%%%%%%%%%%%%%%%%%%%%%%%%%%%%%%%%%%%%%%%%%%%%%%%%%%%%%%%%%%%%%%%
%%%%%    Inhomogeneous linear problem (Remark)    %%%%%%%%%%%%%%%%%%%%%
%%%%%%%%%%%%%%%%%%%%%%%%%%%%%%%%%%%%%%%%%%%%%%%%%%%%%%%%%%%%%%%%%%%%%%%
\begin{remark}[Inhomogeneous linear problem.]
\label{rem-Inhomog}\nopagebreak
$\;$
\begingroup\rm
An {\sl\bfseries explicit formula\/} for calculating the weak solution
$u_{m+1}\colon \RR^1\times [0,T]\to \mathbb{R}$
to problem \eqref{e:BS_u_m+1(t)=f}, \eqref{e:BS_M=u_m+1,0}
will be described later in Section~\ref{s:Numer_parab},
Corollary~\ref{cor-Evolution}.
Numerical methods for computing this solution will be discussed
in {\rm Remark~\ref{rem-MonteCarlo}}, as well.
\null\hfill\Square
%\hfill\Square
\endgroup
\end{remark}
%%%%%%%%%%%%%%%%%%%%%%%%%%%%%%%%%%%%%%%%%%%%%%%%%%%%%%%%%%%%%%%%%%%%%%%
\par\vskip 10pt

As an obvious consequence of our construction we conclude that,
in addition to $u_0$, also each iterate $u_j$; $j = 1,2,3,\dots$,
is a (weak) supersolution of problem
\eqref{e:BS_v(t)}, \eqref{e:BS_M=v,0}.

Standard application of Lebesgue's monotone (or dominated)
convergence theorem yields the monotone pointwise convergence
$u_m(x,t)\searrow v(x,t)$ as $m\nearrow \infty$
for a.e.\ $(x,t)\in \RR^1\times [0,T]$, as well as
the $L^2$\--type norm convergence
$\| u_m(t) - v(t)\|_H\searrow 0$ as $m\nearrow \infty$,
for a.e.\ $t\in (0,T)$, where
$v(t)\equiv v( \,\cdot\, ,t)\in H$ satisfies
\begin{equation}
\label{e:v<=u_m}
\begin{aligned}
  {}- u_0(x,t)\leq v(x,t)\leq u_m(x,t)\leq u_0(x,t)
    \quad\mbox{ for }\, (x,t)\in \RR^1\times [0,T] \,.
\end{aligned}
\end{equation}
Furthermore, we get another $L^2$\--type norm convergence
in the Lebesgue\-(-Hilbert) space\hfil\break
$L^2([0,T]\to H)$, namely,
\begin{equation}
\label{e:L^2:|u_m-v|}
\begin{aligned}
\textstyle
  \| u_m-v\|_{ L^2([0,T]\to H) } \eqdef
    \left( \int_0^T \| u_m(t) - v(t)\|_H^2 \,\mathrm{d}t
    \right)^{1/2} \;\searrow\; 0
    \quad\mbox{ as }\, m\nearrow \infty \,.
\end{aligned}
\end{equation}
We combine this result with Theorem 1.2 (and its proof)
in {\sc A.\ Pazy} \cite[Chapt.~6, {\S}6.1, pp.\ 184--185]{Pazy}
to conclude that $v\in C([0,T]\to H)$ and
$v\colon \RR^1\times [0,T]\to \RR$
is a {\it\bfseries mild solution\/} to the initial value problem
\eqref{e:BS_v(t)}, \eqref{e:BS_M=v,0}.
Finally, applying the well\--known results from
 {\sc L.~C.\ Evans} \cite{Evans-98}, Chapt.~7, {\S}1.1, p.~352, or
 {\sc J.-L.\ Lions} \cite{Lions-61}, Chapt.~IV, {\S}1, p.~44, or
 \cite{Lions-71}, Chapt.~III, eq.~(1.11), p.~102,
we find out that $v$ is a (weak) solution to problem
\eqref{e:BS_v(t)}, \eqref{e:BS_M=v,0}.

We summarize the most important results from our
{\rm monotone iteration scheme\/}
\eqref{e:u_0:Supersol} -- \eqref{e:u_m+1<=u_m}
for problem
\eqref{e:BS_v(t)}, \eqref{e:BS_M=v,0}
in the following theorem.
Precise definitions of the H{\"o}lder spaces used in this theorem,
\begin{math}
  H^{\theta, \theta /2}\bigl( D_{1+1}^{(T)} \bigr)
\end{math}
(a local H{\"o}lder space) and
\begin{math}
  H^{2+\theta, 1 + (\theta /2)}(\overline{Q'})
\end{math}
over the parabolic domain $D_{1+1}^{(T)} = \RR^1\times (0,T)$
(an open strip in $\RR^1\times \RR$)
and its compact subset
$\overline{Q'} = [a,b]\times [\tau, T']\subset D_{1+1}^{(T)}$
(a compact rectangle), respectively, with
$\theta\in (0,1)$,
$-\infty < a < b < +\infty$, and $0 < \tau < T'< T$,
can be found in
{\sc O.~A.\ Ladyzhenskaya}, {\sc V.~A.\ Solonnikov}, and
{\sc N.~N. Ural'tseva}
\cite[Chapt.~I, {\S}1, pp.\ 7--8]{Ladyzh-Ural'-Solonn}.

%%%%%%%%%%%%%%%%%%%%%%%%%%%%%%%%%%%%%%%%%%%%%%%%%%%%%%%%%%%%%%%%%%%%%%%
%%%%%    Main Theorem (Theorem)    %%%%%%%%%%%%%%%%%%%%%%%%%%%%%%%%%%%%
%%%%%%%%%%%%%%%%%%%%%%%%%%%%%%%%%%%%%%%%%%%%%%%%%%%%%%%%%%%%%%%%%%%%%%%
\begin{theorem}[Monotone iterations.]
\label{thm-Iteration}\nopagebreak
$\;$
Let\/ $v_0\in H$ obey\/
{\rm Hypothesis\/} \eqref{hy:v_0} with ineq.~\eqref{e:BS_v(0)}.
Then the monotone iterations\/
\begin{math}
  u_0\geq u_1\geq \dots\geq u_{j-1}\geq u_j\geq \dots\geq {}- u_0 \,,
\end{math}
described in items {\rm (a)}, {\rm (b)}, and {\rm (c)} above,
converge in the Lebesgue\-(-Hilbert) space
$L^2([0,T]\to H)$ to a function
$v\colon \RR^1\times (0,T)$ according to formula
\eqref{e:L^2:|u_m-v|}.
The limit function, $v\in L^2([0,T]\to H)$,
is a (weak) solution to problem
\eqref{e:BS_v(t)}, \eqref{e:BS_M=v,0}.
Furthermore, there is a constant $\theta\in (0,1)$ such that\/
$u_m\in H^{\theta, \theta /2}\bigl( D_{1+1}^{(T)} \bigr)$
holds for every $m = 1,2,3,\dots$, and\/
$v\in H^{\theta, \theta /2}\bigl( D_{1+1}^{(T)} \bigr)$,
as well.

Finally, assume that the function\/
\begin{equation}
%\label{hy:measur}
\nonumber
\tag{\bf G1'}
\begin{aligned}
& G(v;x, \,\cdot\,)\colon t\mapsto G(v;x,t)
                   \colon \RR^1\times (0,T)\to \RR
\\
& \quad\mbox{ is uniformly H{\"o}lder\--continuous on $(0,T)$, }
\\
& \quad\mbox{ uniformly for $(v,x)$ in every bounded subset of }\,
  \RR\times \RR^1 \,.
\end{aligned}
\end{equation}
Then we get even\/
\begin{math}
  u_m\in H^{2+\theta, 1 + (\theta /2)}\bigl( D_{1+1}^{(T)} \bigr)
\end{math}
for every $m = 1,2,3,\dots$, together with
\begin{math}
  v\in H^{2+\theta, 1 + (\theta /2)}\bigl( D_{1+1}^{(T)} \bigr) \,,
\end{math}
where the convergence $u_m\to v$ holds in the norm of the H{\"o}lder space
\begin{math}
  H^{2+\theta', 1 + (\theta'/2)}(\overline{Q'})
\end{math}
over any compact rectangle\/
$\overline{Q'} = [a,b]\times [\tau, T']\subset D_{1+1}^{(T)}$
in the open strip
$D_{1+1}^{(T)} = \RR^1\times (0,T)$, with\/
$-\infty < a < b < +\infty$ and\/ $0 < \tau < T'< T$,
and with any H{\"o}lder exponent $\theta'\in (0,\theta)$.
In particular, each function $u_m$ $(m = 1,2,3,\dots)$
is a {\em strong (classical) solution\/} to problem
\eqref{e:BS_u_j(t)=f}, \eqref{e:BS_M=u_j,0}, whereas $v$
is a {\em strong (classical) solution\/} to problem
\eqref{e:BS_v(t)}, \eqref{e:BS_M=v,0}.
\end{theorem}
%%%%%%%%%%%%%%%%%%%%%%%%%%%%%%%%%%%%%%%%%%%%%%%%%%%%%%%%%%%%%%%%%%%%%%%
%\par\vskip 10pt

%\par\vskip 10pt
%%%%%%%%%%%%%%%%%
\proof
The convergence of our monotone iterations in the Lebesgue\-(-Hilbert) space
\hfil\break
$L^2([0,T]\to H)$ has been established
in eq.~\eqref{e:L^2:|u_m-v|} before this theorem.
Standard regularity theory in
{\sc O.~A.\ Ladyzhenskaya}, {\sc V.~A.\ Solonnikov}, and
{\sc N.~N. Ural'tseva} \cite[Chapt.~III, {\S}10]{Ladyzh-Ural'-Solonn},
Theorem 10.1 on p.~204,
guarantees that there is a constant $\theta\in (0,1)$ such that
$u_m\in H^{\theta, \theta /2}\bigl( D_{1+1}^{(T)} \bigr)$
holds for every $m = 1,2,3,\dots$.
Another regularity result in \cite{Ladyzh-Ural'-Solonn},
Chapt.~IV, {\S}15, Theorem 15.1, pp.\ 405--406,
then yields
\begin{math}
  u_m\in H^{2+\theta, 1 + (\theta /2)}\bigl( D_{1+1}^{(T)} \bigr)
\end{math}
for every $m = 1,2,3,\dots$, with the H{\"o}lder norms
\begin{math}
  \| u_m\|_{ H^{2+\theta, 1 + (\theta /2)}(\overline{Q'}) }
\end{math}
being uniformly bounded, whenever
$\overline{Q'} = [a,b]\times [\tau, T']\subset D_{1+1}^{(T)}$
is a compact rectangle in the open strip
$D_{1+1}^{(T)} = \RR^1\times (0,T)$, with
$-\infty < a < b < +\infty$ and $0 < \tau < T'< T$.
A simple argument based upon the compactness theorem by
{\sc Arzel{\`a}\--Ascoli\/} in
$H^{2+\theta, 1 + (\theta /2)}(\overline{Q'})$
shows that
\begin{math}
  \| v\|_{ H^{2+\theta, 1 + (\theta /2)}(\overline{Q'}) } < +\infty \,,
\end{math}
as well, together with
\begin{math}
  v\in H^{2+\theta, 1 + (\theta /2)}\bigl( D_{1+1}^{(T)} \bigr) \,.
\end{math}
In this context, Theorem 3.1 on pp.\ 983--984 in
{\sc D.~H.\ Sattinger} \cite{Satting-1972}
applies to our 
{\rm monotone iteration scheme\/}
\eqref{e:u_0:Supersol} -- \eqref{e:u_m+1<=u_m}
for problem \eqref{e:BS_v(t)}, \eqref{e:BS_M=v,0}.
%\null\hfill\qed
\qed
%%%%%%%%%%%%%%%%%
\par\vskip 10pt

Two additional important results in
{\sc Sattinger}'s work \cite{Satting-1972},
Theorems 3.3 and 3.4 on p.~986, can be applied to our
{\rm monotone iteration scheme\/},
\eqref{e:u_0:Supersol} -- \eqref{e:u_m+1<=u_m}, as well.
We remark that
the {\sl\bfseries monotonicity (``monotone'') methods\/} in
\cite{Satting-1972}
are formulated for non\-linear elliptic and parabolic
{\sl boundary value problems\/}
in a bounded spatial domain $\Omega\subset \RR^N$.
Nevertheless, they are applicable also to our terminal value problem
\eqref{e:BS_M=hatV}, \eqref{e:BS_M=hatV,T}
(as evidenced by our proof of {\rm Theorem~\ref{thm-Iteration}} above),
even though we have no boundary conditions for the stock price
$S\in (0,\infty)$; more precisely, not for the logarithmic stock price
$x = \log~S\in \RR$ that we use throughout our present article; cf.\
\cite[{\S}3.3, p.~50]{Fencl-Pospis-22}.
In analogy with this {\rm scheme\/},
\eqref{e:u_0:Supersol} -- \eqref{e:u_m+1<=u_m},
for problem \eqref{e:BS_v(t)}, \eqref{e:BS_M=v,0},
which has resulted in the monotone (pointwise) decreasing
sequence of iterations
\begin{math}
  u_0\geq u_1\geq \dots\geq u_{j-1}\geq u_j\geq \dots\geq {}- u_0
\end{math}
a.e.\ in $\RR^1\times (0,T)$ in Theorem~\ref{thm-Iteration} above,
we are able to construct another
{\rm monotone iteration scheme\/}
for problem \eqref{e:BS_v(t)}, \eqref{e:BS_M=v,0}
that is monotone (pointwise) {\sl\bfseries increasing\/} as follows:

%%%%%%%%%%%%%%%%%%%%%%%%%%%%%%%%%%%%%%%%%%%%%%%%%%%%%%%%%%%%%%%%%%%%%%%
%%%%%    Main Theorem (Remark)    %%%%%%%%%%%%%%%%%%%%%%%%%%%%%%%%%%%%%
%%%%%%%%%%%%%%%%%%%%%%%%%%%%%%%%%%%%%%%%%%%%%%%%%%%%%%%%%%%%%%%%%%%%%%%
\begin{remark}[Increasing monotone iterations.]
\label{rem-Iteration}\nopagebreak
$\;$
\begingroup\rm
We set
\begin{equation}
\label{e:w_0:Subsol}
  w_0(x,t) = {}- u_0(x,t)
           = {}- K\, \ee^{\lambda t} \left( \ee^x + \ee^{- x} \right)
    \quad\mbox{ for }\, (x,t)\in \RR^1\times [0,T] \,,
\end{equation}
and recall that $w_0 = {}- u_0$ is a {\em subsolution\/} of problem
\eqref{e:BS_v(t)}, \eqref{e:BS_M=v,0},
by Example~\ref{exam-Sub<=Supersol} with $\kappa = 1$ and some constants\/
$K\geq 1$ and $\lambda\geq 0$ large enough, as specified in this example.
We define the first iteration
$w_1\colon \RR^1\times [0,T]\to \RR$ by replacing the functions
$u_0$ and $u_1$ in eqs.\
\eqref{e:BS_u_1(t)=f} and \eqref{e:BS_M=u_1,0}
by $w_0$ and $w_1$, respectively.
Thus, the inequalities in \eqref{e:u_1<=u_0}
have to be reversed:
\begin{equation}
\label{e:w_1>=w_0}
\begin{aligned}
  {}- u_0(x,t) = w_0(x,t)\leq w_1(x,t)\leq u_0(x,t)
    \quad\mbox{ for }\, (x,t)\in \RR^1\times [0,T] \,.
\end{aligned}
\end{equation}
In order to construct the $j$\--th iteration
$w_j\colon \RR^1\times [0,T]\to \RR$ from
$w_{j-1}\colon \RR^1\times [0,T]\to \RR$; $j = 1,2,3,\dots$,
we replace the pair $(u_{j-1}, u_j)$ by $(w_{j-1}, w_j)$ in 
\eqref{e:u_j<=u_0} -- \eqref{e:u_m+1<=u_m},
thus arriving at the desired monotone (pointwise) increasing
sequence of iterations
\begin{math}
  w_0 = - u_0\leq w_1\leq \dots\leq w_{j-1}\leq w_j\leq \dots\leq u_0
\end{math}
a.e.\ in $\RR^1\times (0,T)$.

The remaing part of {\bf Theorem~\ref{thm-Iteration}\/} remains valid for
$w_m$ in place of $u_m$; $m = 1,2,3,\dots$,
both, {\sl with and without hypothesis\/} {\bf (G1')},
including the monotone convergence $w_m\to v$ in the Lebesgue space
$L^2([0,T]\to H)$ to the function
$v\colon \RR^1\times (0,T)$ described in Theorem~\ref{thm-Iteration}
by formula \eqref{e:L^2:|u_m-v|} and
$u_m\to v$ in the norm of the H{\"o}lder space
\begin{math}
  H^{2+\theta', 1 + (\theta'/2)}(\overline{Q'}) \,,
\end{math}
respectively.
\null\hfill\Square
%\hfill\Square
\endgroup
\end{remark}
%%%%%%%%%%%%%%%%%%%%%%%%%%%%%%%%%%%%%%%%%%%%%%%%%%%%%%%%%%%%%%%%%%%%%%%
\par\vskip 10pt

{\sl Numerical approximation methods\/}
for semi\-linear parabolic systems based on a standard
{\sl finite difference discretization\/} that produces
a {\sl monotone iteration scheme\/}
can be found in the monographs by
{\sc A.~W.\ Leung} \cite{Leung-1989, Leung-2009}.
{\sl Accelerated Monotone Convergence\/}
of this iteration scheme is treated in
\cite{Leung-1989}, Chapt.~VI, {\S}6.3, pp.\ 292--300,
in the context of a two\--point boundary value problem
followed immediately by {\sl $L^2$\--Convergence\/}
for finite\--difference solutions in two dimensional domains
in {\S}6.4, pp.\ 300--323, which can be viewed as
a numerical implementation of our $L^2$\--type norm convergence result
in eq.\ \eqref{e:L^2:|u_m-v|}.
A full, large parabolic system is treated numerically in
\cite[Chapt.~VII]{Leung-1989}, {\S}7.1 -- {\S}7.3, pp.\ 327--357.
Further applications of {\sl monotone iteration schemes\/}
appear in the newer monograph by {\sc A.~W.\ Leung}
\cite{Leung-2009}, Chapt.~3, {\S}3.2 -- {\S}3.3, pp.\ 210--257.

An {\rm important question\/} in this
{\sl\bfseries ``monotone'' approximation method\/}
(cf.\ {\sc David H.\ Sattinger} \cite{Satting-1972})
is the {\sl\bfseries speed of convergence\/}
in eq.~\eqref{e:L^2:|u_m-v|} above.
A suitable norm for estimating this convergence in the Banach space
\begin{math}
  C([0,T]\to H)\hookrightarrow \hfil\break
  L^2((0,T)\to H)
\end{math}
is the {\it\bfseries ``weighted'' supremum norm\/}
\begin{equation}
\label{e:C-weight_norm}
\begin{aligned}
  \vertiii{u}_{ L^{\infty}(0,T) }\eqdef
  \sup_{0\leq t\leq T}
         \left( \ee^{-\omega t}\, \| u(\,\cdot\,,\, t)\|_H \right)
  \hspace*{10pt} ({} < \infty) \,,
\end{aligned}
\end{equation}
where $\omega > 0$ is a sufficiently large positive real number
that guarantees the estimate
\begin{align*}
&   \vertiii{u_{m+1} - u_j}_{ L^{\infty}(0,T) }
  \leq c_{\omega}\,
    \vertiii{u_m - u_{j-1}}_{ L^{\infty}(0,T) }
  \leq \;\dots\;
\\
& \leq c_{\omega}^j\,
    \vertiii{u_{m+1-j} - u_0}_{ L^{\infty}(0,T) }
    \quad\mbox{ for every }\, j = 1,2,3,\dots, m\,;
\end{align*}
$m = 1,2,3,\dots$, with some positive constant $c_{\omega} < 1$.
We refer to
{\sc A.\ Pazy} \cite[Chapt.~2, {\S}2.2, pp.\ 44--45]{Pazy}
for details about how to estimate the constant $\omega > 0$
(from below)
from the values of $r_G = r + L_{\tilde{F}}$ and $L_G$
in eq.~\eqref{e:BS_M=v(t)} and ineq.~\eqref{e:G-Lip}, respectively.
Then also the estimate
\begin{align*}
&   \vertiii{v - u_j}_{ L^{\infty}(0,T) }
  \leq c_{\omega}\,
    \vertiii{v - u_{j-1}}_{ L^{\infty}(0,T) }
  \leq \;\dots\;
\\
& \leq c_{\omega}^j\,
    \vertiii{v - u_0}_{ L^{\infty}(0,T) }
    \quad\mbox{ for every }\, j = 1,2,3,\dots \,,
\end{align*}
follows immediately by letting $m\to \infty$ above.

%%%%%%%%%%%%%%%%%%%%%%%%%%%%%%%%%%%%%%%%%%%%%%%%%%%%%%%%%%%%%%%%%%%%%%%
%%%%%    Numerical mathods: finite diff. & elem. vs. Monte-Carlo    %%%
%%%%%%%%%%%%%%%%%%%%%%%%%%%%%%%%%%%%%%%%%%%%%%%%%%%%%%%%%%%%%%%%%%%%%%%

\section{Numerical methods: finite differences/elements\\
         compared to Monte Carlo}
\label{s:Numer_method}

The Monte Carlo method enjoys broad popularity especially among
the specialists in Probability interested in numerical simulations.
Since its introduction into the toolbox of computational finance
in $1977$ by
{\sc P.\ Boyle\/}, {\sc M.\ Broadie\/}, and {\sc P.\ Glasserman\/}
\cite{Boyle-Broadie-Glass},
{\em\bfseries Monte Carlo\/} simulation has continued to gain acceptance
as a computational method of pricing more and more sophisticated options
and as a risk management tool, as well.
The Monte Carlo method is characterized by its very high flexibility
and its ability to process a multidimensional problem.
One of the main strengths of the Monte Carlo method is that
it does not require discretization, and is therefore
insensitive to space dimension.
Consequently, the more dimensions the problem to be solved has,
the more competitive the Monte Carlo method is when compared to
methods based on space discretization.
However, the rate of convergence of the method,
apparently dependent on the space dimension of the problem
(often referred to as {\sl ``the curse of dimensionality''\/}
 in case of a ``high'' space dimension, which is better known as
 {\sl degeneracy properties of likelihood ratios (LRs)\/}),
is known to be quite low; cf.\
{\sc P.~W.\ Glynn} and {\sc R.~Y.\ Rubinstein} \cite{Glynn-Rubin-2009}.
An important identified limitation of the Monte Carlo method is
its inability
\begingroup\bf
to treat discretization of {\sl non\-linear\/} integral equations
\endgroup
\color{black}
directly.
However, iteration methods that require the solution of a linear problem
at each iteration step (such as ours in Section~\ref{s:Monotone})
have been studied intensively in the past, see e.g.\
{\sc M.~Yu. Plotnikov\/} \cite{Plotnikov-1994}.
The Monte Carlo method is applied there to solve linear problems
in order to construct a {\sl sequence of iterates\/}
that are expected to converge to the solution of
the given {\sl non\-linear\/} integral equation.
This is the strategy we have followed also in the work reported here
(in Section~\ref{s:Monotone}).

Some recent works open perspectives to overcome the difficulties with
discretization of {\sl non\-linear\/} integral equations.
Significant progress has been made on the construction of methods based on
a {\sl\bfseries branching diffusion process\/}.
As the non\-linearities to be treated in numerical discretizations
of PDEs in Finance are mostly {\em\bfseries Lipschitz\--continuous\/}
({\em\bfseries Lipschitzian\/}, for short),
often {\em piecewise linear and continuous\/},
the most efficient {\em\bfseries Monte Carlo\/} algorithms focus
on treating precisely such Lipschitzian non\-linearities.
Among the most interesting earlier pioneering works,
we would like to mention the article by
{\sc P.\ Henry\--Labord{\`e}re\/} in \cite{Labord-2012}.
Here, the author approximates a simple piecewise linear non\-linearity
by a polynomial in order to be able to take advantage of
the probability\--based computational technique,
a {\sl branching diffusion process\/}.
This technique was first introduced by
{\sc H.~P.\ McKean\/} \cite{McKean-1975}
and
{\sc A.~V.\ Skorokhod\/} \cite{Skorokh-1964}
to provide a proba\-bi\-listic representation for the solution of
the {\sl Kolmogorov\--Petrovskii\--Piskunov\/} PDE
({\sl\bfseries KPP\/} equation, for short) and,
more generally, of semi\-linear PDEs whose nonlinearity is
{\sl a power series with non\-negative coefficients\/}
from the interval $[0,1]$.
Since the KPP equation in \cite{McKean-1975, Skorokh-1964}
has {\sl only\/} a quadratic or cubic nonlinearity,
numerical approximation using branching diffusion process
is quite efficient.
However, for applications to
the {\sl\bfseries counterparty risk valuation model\/}
in finance treated in \cite{Labord-2012}, the non\-linearity
that one is interested in is {\sl\bfseries not\/} a polynomial.
It is therefore interesting to investigate methods that can treat
directly {\it\bfseries monotone nonlinearities\/}
that are not necessarily polynomial.
As to the speed and precision of the PDE and Monte\--Carlo methods,
a brief comparison is given in the work by
{\sc G.\ Loeper\/} and {\sc O.\ Pironneau\/} \cite{Loeper-Pironn},
together with a mixed PDE/Monte\--Carlo method that provides better results
for the Heston stochastic volatility model.
An entirely different approach to
{\sc S.~L.\ Heston}'s model \cite{Heston},
based on ``orthogonal series expansion'', is employed in
{\sc B.\ Alziary} and {\sc P.\ Tak\'a\v{c}}
\cite[Sect.~11 (Appendix), {\S}11.1 -- {\S}11.2]{AlziaryTak}, pp.\ 48--50,
in {\sc P.\ Tak\'a\v{c}} \cite{Takac-12},
and in the numerical simulations by
{\sc F.\ Baustian}, {\sc K.\ Filipov{\'a}}, and {\sc J.\ Posp\'{\i}\v{s}il}
\cite{Baust-Pospis-19}.
Replacing a piecewise linear non\-linearity by a polynomial
introduces a significant error into the algorithm in \cite{Labord-2012}.

These branching diffusion processes provide a generalization of
the Feynman\--Kac theory for nonlinear PDEs and lead to
a representation of the solutions in the form of expectation
and thus allow to use Monte Carlo methods for simulations.
A significant progress in the application of branching diffusion process
was made in
{\sc P.\ Henry\--Labord{\`e}re\/}, {\sc X.\ Tan\/}, and {\sc N. Touzi\/}
\cite{Labord-Tan-2014}.
Therefore, more recent studies prefer to work directly with
an arbitrary Lipschitzian non\-linearity
when applying branching diffusion process;
cf.\ {\sc B.\ Bouchard\/}, {\sc X.\ Tan\/}, and {\sc X.\ Warin\/}
\cite{Bouchard-Tan}
and {\sc P.\ Henry\--Labord{\`e}re}, {\sc N.\ Oudjane}, {\sc X.\ Tan},
    {\sc N.\ Touzi}, and {\sc X.\ Warin} \cite{Labord-2019}.
More general (non\--Lipschitzian) non\-linearities are approximated by
(optimal) Lipschitzian non\-linearities on a compact set, if necessary.
In this situation, a standard {\sc Picard\/}'s iteration scheme
(from the proof of the {\sc Picard\--Lindel{\"o}f\/} theorem; see e.g.\
 {\sc Ph.\ Hartman\/} \cite[Chapt.~II, {\S}1, pp.\ 8--10]{Hartman})
can be used to approximate the desired solution $v\in L^2([0,T]\to H)$
to problem \eqref{e:BS_v(t)}, \eqref{e:BS_M=v,0},
i.e., the limit function $v$ in our Theorem~\ref{thm-Iteration}.
Also the approximation error and the rate of convergence
are estimated by classical arguments.

Before switching entirely to a discussion of typically ``analytical''
numerical discretizations, we would like to mention also
some important studies where the two methods, i.e.,
{\em\bfseries Monte Carlo\/} and 
{\em\bfseries finite difference/element\/} algorithms,
are either combined
(\cite{Loeper-Pironn})
or compared against each other in numerical simulations
(\cite[Sect.~4]{Fencl-Pospis-22} and \cite{Loeper-Pironn}).
A suitable combination of these two methods performed on the ``right''
problem can lead to significant improvements of
the ``hybrid'' algorithm by
{\sc G.\ Loeper\/} and {\sc O.\ Pironneau\/} \cite{Loeper-Pironn}
in both, precision and speed.
On the other hand, the work in
{\sc F.\ Baustian}, {\sc M.\ Fencl}, {\sc J.\ Posp\'{\i}\v{s}il}, and
{\sc V.\ {\v{S}}v\'{\i}gler} \cite{Fencl-Pospis-22}
treats the same subject as does our present work, with an
{\sl\bfseries inhomogeneous linear parabolic Cauchy problem\/}
\cite[Eq.\ (6), p.~48]{Fencl-Pospis-22}
closely related to ours in eqs.\
\eqref{e:BS_v(t)=f}, \eqref{e:BS_M=v,0}.
An interesting comparison of several kinds of {\sl Numerical Methods\/}
for the linear initial value problem
\eqref{e:BS_v(t)=f}, \eqref{e:BS_M=v,0}
can be found in that work
\cite[Sect.~4, pp.\ 52--54]{Fencl-Pospis-22}.
In fact, even a quantitative comparison of
{\em\bfseries Monte Carlo\/} and {\em\bfseries finite difference\/}
methods with the exact formula
(for the {\em ``diffusion (or heat) equation''\/})
is provided there \cite[Figs.\ 2--4, p.~53]{Fencl-Pospis-22}.
\color{black}
In \cite[Sect.~3]{Fencl-Pospis-22},
{\rm Monte Carlo\/} and 
{\rm finite difference/element\/} algorithms
are applied separately and then compared against each other
in numerical simulations
\cite[Sect.~4, pp.\ 52--54]{Fencl-Pospis-22}.
In general, the comparison can go ``both ways'';
one has to choose among
{\em\bfseries Monte Carlo\/},
{\em\bfseries finite difference/element\/},
and {\em\bfseries ``hybrid''\/} algorithms according to
the particular problem that should be discretized for
a numerical treatment.

In the previous section (Section~\ref{s:Monotone})
we have presented a method suitable for direct applications of
{\em\bfseries finite difference/element\/} algorithms to
{\sl\bfseries inhomogeneous linear parabolic initial value problems\/}
of type \eqref{e:BS_v(t)=f}, \eqref{e:BS_M=v,0} with
the inhomogeneity (on the right\--hand side in \eqref{e:BS_v(t)=f})
given by the term
$f(x,t) = G(v(x,t); x,t)$ from eq.~\eqref{e:BS_v(t)}.
Our functional\--analytic approach with
an arbitrary Lipschitzian non\-linearity
suggests an analogous {\sc Picard\/} iteration scheme
combined with a finite difference (or similar) discretization method
for the inhomogeneous linear parabolic initial value problem
\eqref{e:BS_v(t)=f}, \eqref{e:BS_M=v,0}.
Solving this problem is especially easy if all functions
$\sigma\colon [0,T]\to (0,\infty)$ and
$q_S, \gamma_S\colon [0,T]\to \RR$
in {\rm Hypotheses\/} \eqref{hy:sigma} and \eqref{hy:q,gamma},
respectively, are constants independent from time $t\in [0,T]$, i.e.,
$\sigma\in (0,\infty)$ and $q_S, \gamma_S\in \RR$.
Our choice of the constant $r_G\in \RR_+$ guarantees that
{\rm Hypothesis\/} \eqref{hy:increasing} is valid:
For almost every pair $(x,t)\in \RR^1\times (0,T)$, the function
\begin{math}
  G(\,\cdot\, ; x,t)\colon v\mapsto G(v; x,t)\colon \RR\to \RR
\end{math}
is monotone increasing.
Under these {\rm Hypotheses\/},
\eqref{hy:sigma} and \eqref{hy:q,gamma},
supplemented by {\rm Hypothesis\/} \eqref{hy:f(t)} in place of
{\rm Hypotheses\/} \eqref{hy:Lipschitz} and \eqref{hy:Hoelder_t}
used for the special inhomogeneity
$f(x,t) = G(v(x,t); x,t)$
on the right\--hand side in eq.~\eqref{e:BS_v(t)},
in the next section (Section~\ref{s:Numer_parab})
we give an {\sl\bfseries explicit formula\/} for the solution $v(x,t)$
to the inhomogeneous linear parabolic initial value problem
\eqref{e:BS_v(t)=f}, \eqref{e:BS_M=v,0}.
This formula is obtained in terms of a linear integral operator
whose kernel is expressed through a ``modification'' of the standard
{\em ``diffusion (or heat) kernel''\/}, cf.\
Theorem~\ref{thm-Evolution} below.

%%%%%%%%%%%%%%%%%%%%%%%%%%%%%%%%%%%%%%%%%%%%%%%%%%%%%%%%%%%%%%%%%%%%%%%
%%%%%    Numerical mathods: semilinear parabolic Cauchy problem    %%%%
%%%%%%%%%%%%%%%%%%%%%%%%%%%%%%%%%%%%%%%%%%%%%%%%%%%%%%%%%%%%%%%%%%%%%%%

\section{Numerical methods: the semilinear parabolic problem}
\label{s:Numer_parab}

In order to calculate (and compute numerically)
the monotone iterations\/
\begin{math}
  u_0\geq u_1\geq \dots\geq u_{j-1}\geq u_j\geq \dots\geq {}- u_0 \,,
\end{math}
treated in Theorem~\ref{thm-Iteration},
we need to solve the inhomogeneous linear parabolic initial value problem
\eqref{e:BS_v(t)}, \eqref{e:BS_M=v,0}
under the assumption that the inhomogeneity
$f(x,t) = G(v(x,t); x,t)$
on the right\--hand side in eq.~\eqref{e:BS_v(t)} is known.
In other words we wish to supplement our main result,
Theorem~\ref{thm-Iteration}, by a result on solving
the initial value problem \eqref{e:BS_v(t)=f}, \eqref{e:BS_M=v,0}.
We will obtain a unique {\it\bfseries weak\/}
(or {\it\bfseries mild\/}) {\it\bfseries solution\/}
$v\colon [0,T]\to H$ under {\rm Hypothesis\/} \eqref{hy:f(t)},
that is to say,
\begin{math}
  f\colon [0,T]\to H
\end{math}
is a continuous function, where
\begin{math}
  f\colon \RR^1\times [0,T]\to \RR\colon (x,t)\mapsto f(x,t)
\end{math}
satisfies
$f(t)\equiv f( \,\cdot\, ,t)\in H$ for every $t\in [0,T]$.

Furthermore, our special inhomogeneity choice of
$f(x,t) = G(v(x,t); x,t)$
satisfies even the stronger hypothesis \eqref{e:f-Hoelder}
on the {\it H{\"o}lder continuity\/} of the function
$f\colon [0,T]\to H$ with the {\it H{\"o}lder exponent\/}
$\vartheta_f\in (0,1)$, see also ineq.~\eqref{e:f(x,t)-Hoelder}.
As a consequence, our {\it\bfseries mild solution\/}
$v\colon [0,T]\to H$, which is continuous by definition,
becomes a unique {\it\bfseries classical solution\/}.
Thus, from now on let us assume that
$f\colon [0,T]\to H$ is a H{\"o}lder\--continuous function
with the {\it H{\"o}lder exponent\/}
$\vartheta_f\in (0,1)$; cf.\ ineq.~\eqref{e:f-Hoelder}
and also ineq.~\eqref{e:f(x,t)-Hoelder} related to
the linear initial value problem
\eqref{e:BS_v(t)=f}, \eqref{e:BS_M=v,0}.

To solve the initial value problem \eqref{e:BS_v(t)=f}, \eqref{e:BS_M=v,0},
we consider the abstract differential equation \eqref{e:abstr:BS_v(t)=f}
with the initial condition \eqref{e:abstr:BS_v_0}, i.e.,
$v(0) = v_0\in H$ at time $t=0$.
A unique {\it\bfseries mild solution\/}
$v\in C([0,T]\to H)$ to this abstract initial value problem
is obtained by the following
{\sl\bfseries variation\--of\--constants\/} formula,
\begin{equation}
\label{e:v(t,0)=exp(I+A_1+A_2)}
\begin{aligned}
  v(t) = \mathfrak{T}(t,0) v_0
       + \int_0^t \mathfrak{T}(t,\tau)\, f(\tau) \,\mathrm{d}\tau
    \quad\mbox{ in $H$ for }\, t\in [0,T] \,,
\end{aligned}
\end{equation}
where
\begin{math}
  \mathfrak{T} =
  \{ \mathfrak{T}(t,s)\colon H\to H\colon 0\leq s\leq t\leq T\}
\end{math}
stands for the {\em\bfseries evolutionary family\/}
of bounded linear operators
\begin{math}
  \mathfrak{T}(t,s)\colon H\to H\ (0\leq s\leq t\leq T)
\end{math}
that yield the unique {\it\bfseries mild solution\/}
\begin{equation*}
%\label{e:v(t,s)=exp(I+A_1+A_2)}
\begin{aligned}
  v(t) = \mathfrak{T}(t,s) v_s
    \quad\mbox{ in $H$ for }\, t\in [s,T]
\end{aligned}
\end{equation*}
to the homogeneous abstract differential equation
(with $f(t)\equiv 0$) on the interval $[s,T]$ for each $s\in [0,T)$,
\begin{alignat}{2}
\label{e:r_G:BS_v(t)=0}
  \frac{\partial v}{\partial t} - \mathcal{A}(t) v
  + r_G\, v
& {}
  = 0
&&  \quad\mbox{ in $H$ for }\, t\in (s,T) \,;
\\
\label{e:abstr:BS_v_s}
    \quad\mbox{ with }\quad
    v(s) & {}= v_s\in H
&&  \quad\mbox{ (the initial condition).}
\end{alignat}
We refer the reader to the monograph by
{\sc A.\ Pazy} \cite[Chapt.~5, {\S}5.6]{Pazy}, Theorem 6.8 on p.~164,
for greater details on the {\em evolutionary family\/}
$\mathfrak{T}$ for the homogeneous abstract initial value problem
\eqref{e:r_G:BS_v(t)=0}, \eqref{e:abstr:BS_v_s}
and to
\cite[Chapt.~5, {\S}5.7]{Pazy}, Theorem 7.1 on p.~168, for the
{\sl variation\--of\--constants\/} formula
\eqref{e:v(t,0)=exp(I+A_1+A_2)}
applied to the in\-homogeneous abstract initial value problem
\eqref{e:abstr:BS_v(t)=f}, \eqref{e:abstr:BS_v_0}
posed on the entire time interval $[0,T]$.

In our theorem below, Theorem~\ref{thm-Evolution},
we provide an explicit formula for
the {\rm evolutionary family\/} $\mathfrak{T}$.
Each operator
\begin{math}
  \mathfrak{T}(t,s)\colon H\to H\ (0\leq s < t\leq T)
\end{math}
turns out to be a (bounded) integral operator with
the {\em\bfseries kernel\/} $\mathfrak{K}(x,y;t,s) > 0$
described explicitly in terms of the standard
{\em ``diffusion (or heat) kernel''\/}
\begin{align}
\label{e:heat_kernel}
\begin{aligned}
  \mathfrak{G}(x;t)\eqdef
  \frac{1}{\sqrt{4\pi t}}\, \exp\left( {}- \frac{x^2}{4t} \right)
  \hspace*{10pt} ({} > 0)
    \quad\mbox{ for }\; (x,t)\in \RR^1\times (0,\infty) \,;
\end{aligned}
\\
\nonumber
    \quad\mbox{ hence, we have }\quad
\textstyle
  \int_{-\infty}^{+\infty} \mathfrak{G}(x;t) \,\mathrm{d}x = 1
    \quad\mbox{ for every }\, t\in (0,\infty) \,.
\end{align}
We refer to
{\sc F.\ John\/} \cite[Chapt.~7, {\S}1(a), pp.\ 206--213]{John}
for greater details about the {\em ``diffusion (or heat) equation''\/}.
Next, we make use of the abbreviations
\begin{equation}
\label{e:S(t),R(t)}
\textstyle
  \mathscr{S}(t)\eqdef \frac{1}{2}\int_0^t [\sigma(\tau)]^2 \,\mathrm{d}\tau
    \quad\mbox{ and }\quad
  \mathscr{R}(t)\eqdef \int_0^t
    \left( q_S(\tau) - \gamma_S(\tau) - \frac{1}{2}\, [\sigma(\tau)]^2
    \right) \,\mathrm{d}\tau
\end{equation}
defined for every $t\in [0,T]$.
By our {\rm Hypothesis\/} \eqref{hy:sigma}, we have
\begin{math}
  \sigma_0 = \min_{t\in [0,T]} \sigma(t) > 0
\end{math}
which entails
\begin{equation*}
\textstyle
    \mathscr{S}(t) - \mathscr{S}(s)
  = \frac{1}{2}\int_s^t [\sigma(\tau)]^2 \,\mathrm{d}\tau
  \geq \frac{1}{2}\, \sigma_0^2\, (t-s) > 0
    \quad\mbox{ whenever }\; 0\leq s < t\leq T \,.
\end{equation*}
%

%%%%%%%%%%%%%%%%%%%%%%%%%%%%%%%%%%%%%%%%%%%%%%%%%%%%%%%%%%%%%%%%%%%%%%%
%%%%%    Inhomogeneous linear problem (Theorem)    %%%%%%%%%%%%%%%%%%%%
%%%%%%%%%%%%%%%%%%%%%%%%%%%%%%%%%%%%%%%%%%%%%%%%%%%%%%%%%%%%%%%%%%%%%%%
\begin{theorem}[Homogeneous linear problem.]
\label{thm-Evolution}\nopagebreak
$\;$
Under the {\rm Hypotheses\/} \eqref{hy:sigma} and\/ \eqref{hy:q,gamma}
stated in {\rm Section~\ref{s:Funct_anal}\/},
the {\em\bfseries evolutionary family\/} $\mathfrak{T}$
on the Hilbert space $H$ takes the following form:
\begin{equation}
\label{e:K(t,s)=exp(I+A_1+A_2)}
\begin{aligned}
& v(t)
  = \left[ \mathfrak{T}(t,s) v_s\right] (x)
  = \int_{-\infty}^{+\infty} \mathfrak{K}(x,y;t,s)\, v_s(y) \,\mathrm{d}y
\\
&   \quad\mbox{ for all $x\in \RR^1$, $0\leq s < t\leq T$, and }\,
    v_s\in H \,,
\end{aligned}
\end{equation}
with the {\em\bfseries kernel\/}
$\mathfrak{K}(x,y;t,s) > 0$ defined by
\begin{align}
\label{e:kernel_exp(I+A_1+A_2)}
\begin{aligned}
&   \mathfrak{K}(x,y;t,s) \equiv \mathfrak{K}(x-y;t,s)
\\
& {}
  = \exp( {}- r_G (t-s) )\cdot
    \mathfrak{G}
    \left( x-y + \mathscr{R}(t) - \mathscr{R}(s) ;\,
                 \mathscr{S}(t) - \mathscr{S}(s)
    \right)
\\
& {}
  = \frac{ \exp( {}- r_G (t-s) ) }%
         { \sqrt{ 4\pi [ \mathscr{S}(t) - \mathscr{S}(s) ] } }\,
    \exp
    \left(
  {}- \frac{( x-y + \mathscr{R}(t) - \mathscr{R}(s) )^2}%
           { 4 [ \mathscr{S}(t) - \mathscr{S}(s) ] }
    \right)
\end{aligned}
\\
\nonumber
    \quad\mbox{ for all $x,y\in \RR^1$ and $0\leq s < t\leq T$. }
\end{align}
Here, the continuous function $v\colon [s,T]\to H$
defined in formula \eqref{e:K(t,s)=exp(I+A_1+A_2)} for $t\in (s,T]$,
with $v(s) = v_s\in H$,
is the unique {\it\bfseries mild solution\/}
of the homogeneous abstract initial value problem
\eqref{e:r_G:BS_v(t)=0}, \eqref{e:abstr:BS_v_s}
on the interval $[s,T]$.
In fact, $v\colon [s,T]\to H$ is also a {\it\bfseries classical solution\/}
to this problem.
\end{theorem}
%%%%%%%%%%%%%%%%%%%%%%%%%%%%%%%%%%%%%%%%%%%%%%%%%%%%%%%%%%%%%%%%%%%%%%%
\par\vskip 10pt

In fact, our existence and uniqueness result in this theorem
corresponds to that in 
{\sc A.\ Pazy} \cite[Chapt.~5, {\S}5.6]{Pazy}, Theorem 6.8 on p.~164.
However, since we wish to calculate the kernel
$\mathfrak{K}(x,y;t,s)$ in formula \eqref{e:kernel_exp(I+A_1+A_2)}
explicitly, we cannot use the abstract proof of
\cite[Chapt.~5, Theorem 6.8]{Pazy} directly.
Thus, we prefer to calculate
the {\rm evolutionary family\/} $\mathfrak{T}$
directly by taking advantage of well\--known results
\cite[Chapt.~7, {\S}1(a), pp.\ 206--213]{John}
for the {\em ``diffusion equation''\/}.

But, before giving the proof of {\bf Theorem~\ref{thm-Evolution}\/},
let us state an important corollary concerning the solvability of
the in\-homogeneous abstract initial value problem
\eqref{e:abstr:BS_v(t)=f}, \eqref{e:abstr:BS_v_0}.

%\par\vskip 10pt
%%%%%%%%%%%%%%%%%%%%%%%%%%%%%%%%%%%%%%%%%%%%%%%%%%%%%%%%%%%%%%%%%%%%%%%
%%%%%    Inhomogeneous linear problem (Corollary)    %%%%%%%%%%%%%%%%%%
%%%%%%%%%%%%%%%%%%%%%%%%%%%%%%%%%%%%%%%%%%%%%%%%%%%%%%%%%%%%%%%%%%%%%%%
\begin{corollary}[Inhomogeneous linear problem.]
\label{cor-Evolution}\nopagebreak
$\;$
Under the {\rm Hypotheses\/} \eqref{hy:sigma} and\/ \eqref{hy:q,gamma}
stated in {\rm Section~\ref{s:Funct_anal}\/}, and\/
$f\in C([0,T]\to H)$, the initial value problem
\eqref{e:abstr:BS_v(t)=f}, \eqref{e:abstr:BS_v_0}
possesses a unique {\it\bfseries mild solution\/}
$v\colon [0,T]\to H$ that is given by the
{\em\bfseries variation\--of\--constants\/} formula
\eqref{e:v(t,0)=exp(I+A_1+A_2)}.
If\/ $f\colon [0,T]\to H$ obeys also {\rm Hypothesis\/} \eqref{hy:f(t)}
stated in {\rm Section~\ref{s:Funct_anal}\/}, then
this solution is also a {\it\bfseries classical solution\/} to problem
\eqref{e:abstr:BS_v(t)=f}, \eqref{e:abstr:BS_v_0}.
\end{corollary}
%%%%%%%%%%%%%%%%%%%%%%%%%%%%%%%%%%%%%%%%%%%%%%%%%%%%%%%%%%%%%%%%%%%%%%%
\par\vskip 10pt

This corollary follows directly from a simple combination of our
Theorem~\ref{thm-Evolution} and a result in
{\sc A.\ Pazy} \cite[Chapt.~5, {\S}5.7]{Pazy}, Theorem 7.1 on p.~168.

%\par\vskip 10pt
%%%%%%%%%%%%%%%%%
{\it Proof of\/} {\bf Theorem~\ref{thm-Evolution}.}$\;$
%\proof
The ``elliptic'' differential operator
\begin{math}
  \mathcal{A}(t)\colon D_{\CC}\subset H_{\CC}\to H_{\CC}
\end{math}
defined in formula \eqref{e:BS_oper} with the (complex) domain
$\mathcal{D}(\mathcal{A}(t)) = D_{\CC}$, i.e.,
\begin{align}
%\label{e:A_1+A_2:BS_oper}
\nonumber
&
\begin{aligned}
& (\mathcal{A}(t) v)(x) =
  (\mathcal{A}_1(t) v)(x) + (\mathcal{A}_2(t) v)(x)
\\
& \equiv
    \frac{1}{2}\, [\sigma(t)]^2 \,\frac{\partial^2 v}{\partial x^2}
  + \left( q_S(t) - \gamma_S(t) - \frac{1}{2}\, [\sigma(t)]^2
    \right) \frac{\partial v}{\partial x}
\end{aligned}
\\
\nonumber
&   \quad\mbox{ for }\; v\in D_{\CC}
    \,\mbox{ and every time }\, t\in (0,T) \,,
\end{align}
is a sum of two {\em\bfseries commuting\/}
differential operators (meaning that their resolvents commute),
respectively,
\begin{math}
  \mathcal{A}_1(t) ,\, \mathcal{A}_2(t)
  \colon D_{\CC}\subset H_{\CC}\to H_{\CC} ,
\end{math}
where
\begin{align*}
%\label{e:A_1:BS_oper}
&
\left\{
\quad
\begin{aligned}
& (\mathcal{A}_1(t) v)(x) \eqdef
    \frac{1}{2}\, [\sigma(t)]^2 \,\frac{\partial^2 v}{\partial x^2}
    \quad\mbox{ for }\; x\in \RR^1
\\
&   \quad\mbox{ is a diffusion operator in $D_{\CC}$\hspace*{5.0mm} and }
\end{aligned}
\right.
\\
%\label{e:A_2:BS_oper}
&
\left\{
\quad
\begin{aligned}
& (\mathcal{A}_2(t) v)(x) \eqdef
    \left( q_S(t) - \gamma_S(t) - \frac{1}{2}\, [\sigma(t)]^2
    \right) \frac{\partial v}{\partial x}
    \quad\mbox{ for }\; x\in \RR^1
\\
&   \quad\mbox{ is a convection (or drift) operator in $D_{\CC}$. }
\end{aligned}
\right.
\end{align*}

In order to construct the sum
$\mathcal{A}(t) = \mathcal{A}_1(t) + \mathcal{A}_2(t)$
with the domain
$D_{\CC}\subset H_{\CC}\to H_{\CC}$ in the Hilbert space
$H_{\CC} = L^2(\mathbb{R};\mathfrak{w})$,
one may apply a rather general perturbation result from
{\sc A.\ Pazy} \cite[Chapt.~3, {\S}3.2]{Pazy}, Theorem 2.1 on p.~80.
Consequently, the (unique) {\em\bfseries classical solution\/}
$v\in C([0,T]\to H)$ to the homogeneous abstract differential equation
\begin{alignat}{2}
\label{e:r_G=0:BS_v(t)=0}
  \frac{\partial v}{\partial t} - \mathcal{A}(t) v
& {}
  = 0
&&  \quad\mbox{ in $H$ for }\, t\in (0,T) \,;
\end{alignat}
with the initial condition $v(0) = v_0\in H$
in eq.~\eqref{e:abstr:BS_v_0},
takes the ``$C^0$\--semigroup'' form
\begin{alignat}{2}
\label{e:v(t)=exp(A_1+A_2)}
&
  v(t) = \mathfrak{T}_{1,2}(t,0) v_0
    \quad\mbox{ in $H$ for }\, t\in [0,T] \,, \quad\mbox{ where }
\\
\label{e:exp(A_1+A_2)}
&
\begin{aligned}
    \mathfrak{T}_{1,2}(t,s)
  = \exp\left(
        [ \mathscr{S}(t) - \mathscr{S}(s) ]\, \frac{\partial^2}{\partial x^2}
        \right)
  \,\exp\left(
        [ \mathscr{R}(t) - \mathscr{R}(s) ]\, \frac{\partial}{\partial x}
        \right)
\\
    \quad\mbox{ in $H$, for }\, 0\leq s\leq t\in [0,T] \,,
\end{aligned}
\end{alignat}
is the {\rm evolutionary family\/}
\begin{math}
  \mathfrak{T}_{1,2} =
  \{ \mathfrak{T}_{1,2}(t,s)\colon H\to H\colon 0\leq s\leq t\leq T\}
\end{math}
of bounded linear operators
\begin{math}
  \mathfrak{T}_{1,2}(t,s)\colon H\to H\ (0\leq s\leq t\leq T)
\end{math}
that yield the unique {\it\bfseries classical solution\/}
$v(t) = \mathfrak{T}_{1,2}(t,s) v_s$ in $H$
to the homogeneous abstract differential equation
\eqref{e:r_G=0:BS_v(t)=0}
(with $f(t)\equiv 0$) on the interval $[s,T]$ for each $s\in [0,T)$,
with the initial condition $v(s) = v_s\in H$.
We recall that the auxiliary functions
\begin{math}
  t \,\mapsto\, \mathscr{S}(t) ,\, \mathscr{R}(t)\colon [0,T]\to \RR
\end{math}
have been defined in eq.~\eqref{e:S(t),R(t)}.
We refer the reader to the general results in
{\sc A.\ Pazy} \cite[Chapt.~5]{Pazy},
Theorem 6.8 in {\S}5.6 on p.~164 and
Theorem 7.1 in {\S}5.7 on p.~168.
It is easy to see that
\begin{align}
\label{e:v'(t)=exp(A_1)}
&
\left\{\quad
\begin{aligned}
& \frac{\mathrm{d}}{\mathrm{d}t}\,
    \exp\left( \mathscr{S}(t)\, \frac{\partial^2}{\partial x^2} \right)
    v_0
  = \mathscr{S}'(t)
    \exp\left( \mathscr{S}(t)\, \frac{\partial^2}{\partial x^2} \right)
    \frac{\mathrm{d}^2 v_0}{\mathrm{d}x^2}
\\
& {}
  = \mathscr{S}'(t)\cdot
    \exp\left( \mathscr{S}(t)\, \frac{\partial^2}{\partial x^2} \right)
    v_0''
  = \frac{1}{2}\, [\sigma(t)]^2\cdot
    \exp\left( \mathscr{S}(t)\, \frac{\partial^2}{\partial x^2} \right)
    v_0''
\\
&   \quad\mbox{ holds for all }\, v_0\in D_{\CC}
    \,\mbox{ and for every }\, t\in (0,T] \,,
    \quad\mbox{ and }\quad
\end{aligned}
\right.
\\
\label{e:v(t)=exp(A_2)}
&
\begin{aligned}
& \left[
    \exp\left( \mathscr{R}(t)\, \frac{\partial}{\partial x} \right)
    v_0
  \right](x) = v_0( x + \mathscr{R}(t) )
    \quad\mbox{ for all $x\in \RR^1$ and $t\in [0,T]$, }
\end{aligned}
\\
\nonumber
&   \qquad\mbox{ together with }\qquad
\\
\label{e:v'(t)=exp(A_2)}
&
\left\{\quad
\begin{aligned}
& \frac{\mathrm{d}}{\mathrm{d}t}\,
    \exp\left( \mathscr{R}(t)\, \frac{\partial}{\partial x} \right)
    v_0
  = \mathscr{R}'(t)
    \exp\left( \mathscr{R}(t)\, \frac{\partial}{\partial x} \right)
    \frac{\mathrm{d} v_0}{\mathrm{d}x}
\\
& {}
  = \mathscr{R}'(t)\cdot v_0'( x + \mathscr{R}(t) )
  = \left( q_S(t) - \gamma_S(t) - \frac{1}{2}\, [\sigma(t)]^2
    \right)\cdot v_0'( x + \mathscr{R}(t) )
\\
&   \quad\mbox{ for all }\, v_0\in D_{\CC}
    \,\mbox{ and for every }\, t\in [0,T] \,.
\end{aligned}
\right.
\end{align}

Finally, to construct the full linear operator
\begin{equation*}
  \mathcal{A}(t) - r_G\,\bullet \colon D_{\CC}\subset H_{\CC}\to H_{\CC}
  \colon v \,\longmapsto\, \mathcal{A}(t) v - r_G\, v
\end{equation*}
that appears on the left\--hand side of
the in\-homogeneous abstract differential equation
\eqref{e:abstr:BS_v(t)=f},
let us combine formula \eqref{e:v(t)=exp(A_1+A_2)} with eqs.\
\eqref{e:v'(t)=exp(A_1)} and \eqref{e:v'(t)=exp(A_2)}
in order to define the {\rm evolutionary family\/}
\begin{math}
  \mathfrak{T} =
  \{ \mathfrak{T}(t,s)\colon H\to H\colon 0\leq s\leq t\leq T\}
\end{math}
by the following composition of bounded linear operators
\begin{align}
\label{e:exp(I+A_1+A_2)}
\begin{aligned}
& \mathfrak{T}(t,s) v_s \eqdef
    \exp( {}- r_G\, (t-s) )\cdot \mathfrak{T}_{1,2}(t,s) v_s
\\
& {}
  = \exp( {}- r_G\, (t-s) )\cdot
    \exp\left(
    [ \mathscr{S}(t) - \mathscr{S}(s) ]\, \frac{\partial^2}{\partial x^2}
        \right)
  \,\exp\left(
    [ \mathscr{R}(t) - \mathscr{R}(s) ]\, \frac{\partial}{\partial x}
        \right) v_s
\end{aligned}
\\
\nonumber
    \quad\mbox{ in $H$ for $v_s\in H$ and }\, 0\leq s\leq t\leq T \,.
\end{align}
Applying eqs.\
\eqref{e:v'(t)=exp(A_1)} and \eqref{e:v'(t)=exp(A_2)}
to this formula we obtain
\begin{equation}
\label{e:(d/dt)exp(I+A_1+A_2)}
\begin{aligned}
&   \frac{\mathrm{d}}{\mathrm{d}t}\, \mathfrak{T}(t,s) v_s
  = \left( \mathcal{A}(t) - r_G\,\bullet \right) \mathfrak{T}(t,s) v_s
\\
&   \quad\mbox{ in $H_{\CC}$ for $v_s\in D_{\CC}$ and }\,
    0\leq s < t\leq T \,.
\end{aligned}
\end{equation}
Let us recall that $f\colon [0,T]\to H$ is a continuous function.
Given any initial value $v_s\in H$, we may apply
\cite[Chapt.~5, {\S}5.7]{Pazy}, Theorem 7.1 on p.~168,
to conclude that the function
\begin{equation}
\label{e:v(t,s)=exp(I+A_1+A_2)}
\begin{aligned}
  v(t) = \mathfrak{T}(t,s) v_s
       + \int_s^t \mathfrak{T}(t,\tau)\, f(\tau) \,\mathrm{d}\tau
    \quad\mbox{ in $H$ for }\, t\in [s,T]
\end{aligned}
\end{equation}
is the unique {\it\bfseries mild solution\/}
of the in\-homogeneous abstract differential equation
\eqref{e:abstr:BS_v(t)=f} in the time interval $[s,T]$
with the initial value $v(s) = v_s\in H$.

In our proof of {\rm Theorem~\ref{thm-Iteration}} above,
we have taken advantage of {\rm Hypotheses\/}
\eqref{hy:Lipschitz} and {\rm ({\bf G1'})}
in order to conclude that the limit function
$v\in L^2([0,T]\to H)$ obtained in formula \eqref{e:L^2:|u_m-v|}
is a {\it\bfseries strong (classical) solution\/} to problem
\eqref{e:BS_v(t)}, \eqref{e:BS_M=v,0}.
Indeed, these two {\rm Hypotheses\/},
\eqref{hy:Lipschitz} and {\rm ({\bf G1'})},
guarantee that all iterates (functions)
$u_j\colon [0,T]\to H$; $j=0,1,2,\dots$,
in {\rm Theorem~\ref{thm-Iteration}}
are uniformly H{\"o}lder\--continuous functions
with some {\it H{\"o}lder exponent\/} $\vartheta_v\in (0,1)$
and their monotone sequence
\begin{math}
  u_0\geq u_1\geq \dots\geq u_{j-1}\geq u_j\geq \dots\geq {}- u_0
\end{math}
is uniformly bounded in the H{\"o}lder space
$C^{\vartheta_v}( [0,T]\to H)$.
Hence, we have
$v\in C^{\vartheta_v}([0,T]\to H)$, as well.
We conclude that our function
$f(x,t) = G(v(x,t); x,t)$ from eq.~\eqref{e:BS_v(t)}
obeys {\rm Hypothesis\/} \eqref{hy:f-Hoelder}.
Consequently, from now on let us assume that
$f\colon [0,T]\to H$ is a H{\"o}lder\--continuous function
with the {\it H{\"o}lder exponent\/}
$\vartheta_f\in (0,1)$; cf.\ ineq.~\eqref{e:f-Hoelder}
and also ineq.~\eqref{e:f(x,t)-Hoelder} related to
the the linear initial value problem
\eqref{e:BS_v(t)=f}, \eqref{e:BS_M=v,0}.
Then, given any initial value $v_s\in H$, it follows again from
\cite[Chapt.~5, {\S}5.7]{Pazy}, Theorem 7.1 on p.~168,
that that the function
$v\colon [s,T]\to H$ defined above by the 
{\sl variation\--of\--constants\/} formula
in \eqref{e:v(t,s)=exp(I+A_1+A_2)},
is even a {\it\bfseries classical solution\/}
of the in\-homogeneous abstract differential equation
\eqref{e:abstr:BS_v(t)=f} in the time interval $[s,T]$
with the initial value $v(s) = v_s\in H$.
From formulas
\eqref{e:v(t)=exp(A_1+A_2)}, \eqref{e:v(t)=exp(A_2)}, and
\eqref{e:v(t,s)=exp(I+A_1+A_2)} we deduce that each mapping
\begin{math}
  \mathfrak{T}(t,s)\colon H\to H\ (0\leq s < t\leq T)
\end{math}
is an integral operator with the {\em\bfseries kernel\/}
$\mathfrak{K}(x,y;t,s) > 0$;
see eq.~\eqref{e:K(t,s)=exp(I+A_1+A_2)}.

Next, we will derive an explicit formula for the kernel
$\mathfrak{K}(x,y;t,s)$ with a help from the standard
{\em ``diffusion kernel''\/}
$\mathfrak{G}(x;t)$ defined in eq.~\eqref{e:heat_kernel}.
First, the standard {\em ``diffusion semigroup''\/}
\begin{math}
  \mathcal{T}_1 =
  \{ \mathcal{T}_1(t)\colon H\to H\colon 0\leq t < \infty\}
\end{math}
consists of (bounded) integral operators
\begin{math}
  \mathcal{T}_1(t) = \exp\left( t\, \frac{\partial^2}{\partial x^2} \right)
  \colon H\to H\ (0\leq t < \infty)
\end{math}
with the kernel $\mathfrak{G}(x-y;t) > 0$, i.e.,
\begin{equation*}
%\label{e:G(t,s)=exp(t.A_1)}
\begin{aligned}
&   \left[ \mathcal{T}_1(t) v_0\right] (x)
  = \int_{-\infty}^{+\infty} \mathfrak{G}(x-y;t)\, v_0(y) \,\mathrm{d}y
\\
&   \quad\mbox{ defined for all $x\in \RR^1$, $0 < t < \infty$, and }\,
    v_0\in H \,.
\end{aligned}
\end{equation*}
We refer to
{\sc F.\ John\/} \cite[Chapt.~7, {\S}1(a), pp.\ 206--213]{John}
for greater details.
Hence, in formula \eqref{e:exp(I+A_1+A_2)}
we have the (bounded) integral operator
\begin{align}
%\label{e:exp(A_1)}
\nonumber
\begin{aligned}
    \left[ \mathfrak{T}_1(t,s) v_s\right] (x)
& {}
  \eqdef
    \left[
    \exp\left(
    [ \mathscr{S}(t) - \mathscr{S}(s) ]\, \frac{\partial^2}{\partial x^2}
        \right) v_s
    \right] (x)
\\
& {}
  = \int_{-\infty}^{+\infty} \mathfrak{G}_1(x-y;t,s)\, v_s(y) \,\mathrm{d}y
\end{aligned}
\\
\nonumber
    \quad\mbox{ defined for all $x\in \RR^1$, $0\leq s < t\leq T$, and }\,
    v_s\in H \,,
\end{align}
with the kernel $\mathfrak{G}_1(x-y;t,s) > 0$ given by
\begin{align}
%\label{e:kernel_exp(A_1)}
\nonumber
\begin{aligned}
&   \mathfrak{G}_1(x-y;t,s)
  = \mathfrak{G}\left( x-y ;\, \mathscr{S}(t) - \mathscr{S}(s) \right)
\\
& {}
  = \frac{1}{ \sqrt{ 4\pi [ \mathscr{S}(t) - \mathscr{S}(s) ] } }\,
    \exp
    \left( {}- \frac{(x-y)^2}{ 4 [ \mathscr{S}(t) - \mathscr{S}(s) ] }
    \right)
\end{aligned}
\\
\nonumber
    \quad\mbox{ for all $x,y\in \RR^1$ and $0\leq s < t\leq T$ }\,.
\end{align}
Let us recall that, by eq.~\eqref{e:S(t),R(t)}, we have
\begin{math}
    \mathscr{S}(t) - \mathscr{S}(s)
  = \frac{1}{2}\int_s^t [\sigma(\tau)]^2 \,\mathrm{d}\tau > 0
\end{math}
whenever $0\leq s < t\leq T$.

Second, the {\em ``shift (or translation) group''\/}
\begin{math}
  \mathcal{T}_2 =
  \{ \mathcal{T}_2(t)\colon H\to H\colon t\in \RR\}
\end{math}
consists of (bounded) translation operators
\begin{math}
  \mathcal{T}_2(t) = \exp\left( t\, \frac{\partial}{\partial x} \right)
  \colon H\to H\ ( -\infty < t < +\infty )
\end{math}
given by the formula
\begin{equation*}
%\label{e:G(t,s)=exp(t.A_2)}
  \left[ \mathcal{T}_2(t) v_0\right] (x) = v_0(x+t)
    \quad\mbox{ for all $x\in \RR^1$, $t\in \RR$, and }\, v_0\in H \,.
\end{equation*}
Hence, on the right\--hand side in formula \eqref{e:exp(I+A_1+A_2)}
we have the translation operator
\begin{align}
%\label{e:exp(A_2)}
\nonumber
\begin{aligned}
    \left[ \mathfrak{T}_2(t,s) v_s\right] (x)
& {}
  \eqdef
    \left[
    \exp\left(
    [ \mathscr{R}(t) - \mathscr{R}(s) ]\, \frac{\partial}{\partial x}
        \right) v_s
    \right] (x)
\\
& {}
  = v_s( x + \mathscr{R}(t) - \mathscr{R}(s) )
\end{aligned}
\\
\nonumber
    \quad\mbox{ for all $x\in \RR^1$, $0\leq s\leq t\leq T$, and }\,
    v_s\in H \,.
\end{align}
Consequently, the (bounded linear) composition operator on $H$,
defined by formula \eqref{e:exp(A_1+A_2)}, i.e.,
\begin{equation*}
%\label{e:exp(A_1+A_2)}
\begin{aligned}
&   \mathfrak{T}_{1,2}(t,s) \eqdef
    \mathfrak{T}_1(t,s)\, \mathfrak{T}_2(t,s) \equiv
    \mathfrak{T}_1(t,s)\circ \mathfrak{T}_2(t,s)
\\
& {}
  = \exp\left(
    [ \mathscr{S}(t) - \mathscr{S}(s) ]\, \frac{\partial^2}{\partial x^2}
        \right)
  \,\exp\left(
    [ \mathscr{R}(t) - \mathscr{R}(s) ]\, \frac{\partial}{\partial x}
        \right) \colon H\to H \,,
\end{aligned}
\end{equation*}
for $0\leq s < t\leq T$, is an integral operator
\begin{align*}
%\label{e:exp(A_1+A_2)}
\begin{aligned}
    \left[ \mathfrak{T}_{1,2}(t,s) v_s\right] (x)
& {}
  = \int_{-\infty}^{+\infty}
    \mathfrak{G}_{1,2}(x,y;t,s)\, v_s(y) \,\mathrm{d}y
\end{aligned}
\\
    \quad\mbox{ defined for all $x\in \RR^1$, $0\leq s < t\leq T$, and }\,
    v_s\in H \,,
\end{align*}
with the kernel $\mathfrak{G}_{1,2}(x,y;t,s) > 0$ given by
\begin{align}
%\label{e:kernel_exp(A_1+A_2)}
\nonumber
\begin{aligned}
&   \mathfrak{G}_{1,2}(x,y;t,s) \equiv \mathfrak{G}_{1,2}(x-y;t,s)
  = \mathfrak{G}_1
    \bigl( x + \mathscr{R}(t) - \mathscr{R}(s) - y ;\, t,s \bigr)
\\
& {}
  = \mathfrak{G}
    \bigl( x-y + \mathscr{R}(t) - \mathscr{R}(s) ;\,
            \mathscr{S}(t) - \mathscr{S}(s)
    \bigr)
\\
& {}
  = \frac{1}{ \sqrt{ 4\pi [ \mathscr{S}(t) - \mathscr{S}(s) ] } }\,
    \exp
    \left(
  {}- \frac{( x-y + \mathscr{R}(t) - \mathscr{R}(s) )^2}%
           { 4 [ \mathscr{S}(t) - \mathscr{S}(s) ] }
    \right)
\end{aligned}
\\
\nonumber
    \quad\mbox{ for all $x,y\in \RR^1$ and $0\leq s < t\leq T$. }
\end{align}

Finally, comparing
\begin{math}
  \mathcal{A}(t)\colon D_{\CC}\subset H_{\CC}\to H_{\CC}
\end{math}
with the full linear operator
\begin{equation*}
  \mathcal{A}(t) - r_G\,\bullet \colon D_{\CC}\subset H_{\CC}\to H_{\CC}
  \colon v \,\longmapsto\, \mathcal{A}(t) v - r_G\, v \,,
\end{equation*}
we arrive at
\begin{math}
  \mathfrak{T}(t,s)
  = \exp( {}- r_G\, (t-s) )\cdot \mathfrak{T}_{1,2}(t,s)
\end{math}
in $H$ for all $0\leq s < t\leq T$.
Hence, eq.~\eqref{e:exp(I+A_1+A_2)} is valid and the desired formulas
\eqref{e:K(t,s)=exp(I+A_1+A_2)} and
\eqref{e:kernel_exp(I+A_1+A_2)} follow.

Our proof of Theorem~\ref{thm-Evolution} is finished.
%\null\hfill\qed
\qed
%%%%%%%%%%%%%%%%%
%\par\vskip 10pt

%\par\vskip 10pt
%%%%%%%%%%%%%%%%%%%%%%%%%%%%%%%%%%%%%%%%%%%%%%%%%%%%%%%%%%%%%%%%%%%%%%%
%%%%%    Inhomogeneous linear problem (Remark)    %%%%%%%%%%%%%%%%%%%%%
%%%%%%%%%%%%%%%%%%%%%%%%%%%%%%%%%%%%%%%%%%%%%%%%%%%%%%%%%%%%%%%%%%%%%%%
\begin{remark}[An alternative proof of Theorem~\ref{thm-Evolution}.]
\label{rem-Evolution}\nopagebreak
$\;$
\begingroup\rm
In order to derive formula \eqref{e:kernel_exp(I+A_1+A_2)}
for the {\em\bfseries kernel\/}
$\mathfrak{K}(x,y;t,s) > 0$ of the integral operator
$\mathfrak{T}(t,s)\colon H\to H$ in
eq.~\eqref{e:K(t,s)=exp(I+A_1+A_2)}, one can replace the calculations
in our proof of Theorem~\ref{thm-Evolution} above
by applying first the Fourier transformation to
the homogeneous parabolic problem
\eqref{e:r_G:BS_v(t)=0}, \eqref{e:abstr:BS_v_s}
to calculate the Fourier transform of the kernel
$\mathfrak{K}(x,y;t,s)\equiv \mathfrak{K}(x-y;t,s)$
with respect to the variable (difference) $x-y\in \RR^1$,
followed by a standard application of
the {\em inverse\/} Fourier transformation to obtain the desired kernel
$\mathfrak{K}(x,y;t,s)$ for all $x,y\in \RR^1$ and $0\leq s < t < \infty$.
This procedure is analogous with the derivation of
the {\em ``diffusion kernel''\/}
$\mathfrak{G}(x;t)$ (defined here in eq.~\eqref{e:heat_kernel})
by Fourier transformation in
{\sc F.~John\/} \cite[Chapt.~7, {\S}1(a), pp.\ 208--209]{John};
see Eq.~(1.10d) on p.~209.
\hfill\Square
\endgroup
\end{remark}
%%%%%%%%%%%%%%%%%%%%%%%%%%%%%%%%%%%%%%%%%%%%%%%%%%%%%%%%%%%%%%%%%%%%%%%
\par\vskip 10pt

Our explicit 
{\sl\bfseries variation\--of\--constants\/} formula
\eqref{e:v(t,0)=exp(I+A_1+A_2)} verified in
{\bf Corollary~\ref{cor-Evolution}\/} (to Theorem~\ref{thm-Evolution})
renders a precise analytic way of calculating the exact solution to problem
\eqref{e:abstr:BS_v(t)=f}, \eqref{e:abstr:BS_v_0},
which is the abstract form of the Cauchy problem
\eqref{e:BS_v(t)=f}, \eqref{e:BS_M=v,0}.
In practical, real world applications to problems in Mathematical Finance
this solution needs to be approximated by suitable Numerical Methods;
typically with reasonable speed and within reasonable precision.
Whereas the {\sl\bfseries non\-linear effects\/},
modelled by the non\-linear ``reaction'' function
\begin{math}
  G(\,\cdot\, ; x,t)\colon v\mapsto G(v; x,t)\colon \RR\to \RR \,,
\end{math}
defined in eq.~\eqref{e:G=F-L_tildeF},
are computed in a most standard way through the in\-homogeneity
$f(x,t) = G(v(x,t); x,t)$ from eq.~\eqref{e:BS_v(t)},
with the (non\-linear) substitution operator
$G(\,\cdot\, ; x,t)\colon \RR\to \RR$,
the numerical computation of the (unique) solution 
$v\colon \RR^1\times (0,T)\to \mathbb{R}$
to the {\sl\bfseries linear\/} initial value Cauchy problem
\eqref{e:BS_v(t)=f}, \eqref{e:BS_M=v,0}
is much more complicated:

%\par\vskip 10pt
%%%%%%%%%%%%%%%%%%%%%%%%%%%%%%%%%%%%%%%%%%%%%%%%%%%%%%%%%%%%%%%%%%%%%%%
%%%%%    Monte Carlo v. finite difference/element (Remark)    %%%%%%%%%
%%%%%%%%%%%%%%%%%%%%%%%%%%%%%%%%%%%%%%%%%%%%%%%%%%%%%%%%%%%%%%%%%%%%%%%
\begin{remark}
% [An alternative proof of Theorem~\ref{thm-Evolution}.]
\label{rem-MonteCarlo}\nopagebreak
\begingroup\rm
$\,$
{\bf (i)}$\,$
Eq.~\eqref{e:BS_v(t)=f} being a linear
{\em ``diffusion equation''\/},
the Monte Carlo method is a very natural way for approximating
the precise analytic solution to problem
\eqref{e:BS_v(t)=f}, \eqref{e:BS_M=v,0}
(i.e., \eqref{e:abstr:BS_v(t)=f}, \eqref{e:abstr:BS_v_0})
by numerical simulations produced by this method.
We refer to the works by
{\sc I.\ Arregui\/}, {\sc B.\ Salvador\/}, and {\sc C.\ V{\'a}zquez\/}
\cite[Sects.\ 3--4, pp.\ 18--23]{Arregui-Salvad},
{\sc F.\ Baustian}, {\sc M.\ Fencl}, {\sc J.\ Posp\'{\i}\v{s}il}, and
{\sc V.\ {\v{S}}v\'{\i}gler}
\cite[{\S}3.6, p.~52]{Fencl-Pospis-22} and to
{\sc M.~Yu. Plotnikov\/} \cite[Sect.~1, pp.\ 121--125]{Plotnikov-1994}
for greater details.
The first two references, \cite{Arregui-Salvad, Fencl-Pospis-22},
treat exactly the problem of
{\it ``Option pricing under some Value Adjustment''\/}
({\sl\bfseries xVA})
in Mathematical Finance; under
{\it ``Credit Value Adjustment''\/}
({\sl\bfseries CVA}), for instance.
The third one, \cite{Plotnikov-1994},
treats linear integral equations of type
\eqref{e:v(t,s)=exp(I+A_1+A_2)}
which envolve the {\rm evolutionary family\/}
of bounded linear integral operators
\begin{math}
  \mathfrak{T}(t,s)\colon H\to H\ (0\leq s < t\leq T)
\end{math}
on the Hilbert space $H$.
In fact, a numerical method for solving the full,
non\-linear integral equation \eqref{e:v(t,0)=exp(I+A_1+A_2)}
for the unknown function $v\in C([0,T]\to H)$
substituted into the (subsequently unknown) non\-linearity
\begin{equation*}
  f(\tau)\equiv f( \,\cdot\, ,\tau)\colon \RR^1\to \RR\colon
  x \,\mapsto\, f(x,\tau) = G(v(x,\tau); x,\tau)
    \quad\mbox{ with }\; f(\tau)\in H
\end{equation*}
for every $\tau\in [0,T]$,
is provided in \cite{Plotnikov-1994}.
This work is based in an iteration method for
a non\-linear integral equation similar to ours, cf.\
\cite[Eq.\ (1.1), p.~121]{Plotnikov-1994}.

{\bf (ii)}$\,$
In contrast with the probabilistic Monte Carlo methods for solving
the Cauchy problem
\eqref{e:BS_v(t)=f}, \eqref{e:BS_M=v,0},
analytic methods based on
a {\em\bfseries finite difference (or finite element) scheme\/}
provide a highly competitive alternative to Monte Carlo
in a series of works, such as the monograph by
{\sc Y.\ Achdou} and {\sc O.\ Pironneau} \cite{Achdou-Pironn}
and the articles by
{\sc I.\ Arregui\/} et {\sc al.}
\cite[{\S}3.4, pp.\ 20--21]{Arregui-Salvad},
\cite[Sect.~4, pp.\ 734--737]{Arregui-Sevco},
{\sc F.\ Baustian} et {\sc al.}
\cite[{\S}3.5, pp.\ 50--51]{Fencl-Pospis-22},
{\sc M.~N.\ Koleva\/} \cite[Sect.~3, pp.\ 367--368]{Koleva-2013}, and
{\sc M.~N.\ Koleva\/} and {\sc L.~G.\ Vulkov}
% \cite[Sect.~2, pp.\ 568--570]{KolevaVul-2012},
% \cite[Sects.\ 2--3, pp.\ 432--434]{KolevaVul-2017}, and
\cite[Sects.\ 3--4, pp.\ 510--515]{KolevaVul-2018}.

{\bf (iii)}$\,$
Last but not least, a ``hybrid'' algorithm mixing
{\em\bfseries Monte Carlo\/} with
{\em\bfseries finite difference/element\/} methods
in quest for optimization on both, precision and speed, is presented in
{\sc G.\ Loeper\/} and {\sc O.\ Pironneau\/} \cite{Loeper-Pironn}.
\hfill\Square
\endgroup
\end{remark}
%%%%%%%%%%%%%%%%%%%%%%%%%%%%%%%%%%%%%%%%%%%%%%%%%%%%%%%%%%%%%%%%%%%%%%%
%\par\vskip 10pt

%%%%%%%%%%%%%%%%%%%%%%%%%%%%%%%%%%%%%%%%%%%%%%%%%%%%%%%%%%%%%%%%%%%%%%%

%%%%%%%%%%%%%%%%%%%%%%%%%%%%%%%%%%%
% \section*{Acknowledgment.}
% A part of this research was performed while
% the second author (P.T.) was a visiting professor at
% Toulouse School of Economics, I.M.T.,
% Universit\'e de Toulouse -- Capitole, Toulouse, France.
%%%%%%%%%%%%%%%%%%%%%%%%%%%%%%%%%%%

%The authors express their thanks to
%an anonymous referee for
%very careful reading of the manuscript and
%a number of helpful suggestions.

%\newpage
%%%%%%%%%%%%%%%%%%%%%%%%%%%%%%%%%%%%%%%%%%%%%%%%%%%%%%%%%%%%%%%%%%%%%%%
%%%%%    BIBLIOGRAPHY    %%%%%%%%%%%%%%%%%%%%%%%%%%%%%%%%%%%%%%%%%%%%%%
%%%%%%%%%%%%%%%%%%%%%%%%%%%%%%%%%%%%%%%%%%%%%%%%%%%%%%%%%%%%%%%%%%%%%%%

%%%%%%%%%%%%%%%%%%%%%%%%%%%%%%%%%%%%%%%%%%%%%%%%%%%%%%%%%%%%%%%%%%%%%%%%%%%%
%               AMS-LaTeX Version 1.1 file for electronic submission       %
% PAPER.BBL                  			          October 1999     %
%                                                                          %
%       Ground\--State Positivity, Negativity, and Compactness		   %
%       for a Schr\"odinger Operator in $\mathbb{R}^N$			   %
%                                                                          %
%%%%%%%%%%%%%%%%%%%%%%%%%%%%%%%%%%%%%%%%%%%%%%%%%%%%%%%%%%%%%%%%%%%%%%%%%%%%
%
%\bibliographystyle{amsplain}        %% if within the text
%
\makeatletter \renewcommand{\@biblabel}[1]{\hfill#1.} \makeatother
\ifx\undefined\bysame
\newcommand{\bysame}{\leavevmode\hbox to3em{\hrulefill}\,}
\fi
%

%
%%%%%%%%%%%%%%%%%%%%%%%%%%%%%%%%%%%%%%%%%%%%%%%%%%%%%%%%%%%%%%%%%%%%%%%
%%%%%%%%%%%%%%%%%%%%%%%%%%%%%%%%%%%%%%%%%%%%%%%%%%%%%%%%%%%%%%%%%%%%%%%
%%%%%%%%%%%%%%%%%%%%%%%%%%%%%%%%%%%%%%%%%%%%%%%%%%%%%%%%%%%%%%%%%%%%%%%
%
%%%%%%%%%%%%%%%%%%%%%%%%%%%
%
\end{document}